\documentclass[11pt,reqno]{amsart}
\usepackage{amssymb}
\usepackage[curve,matrix,arrow]{xy}
\newcommand{\A} {\mathcal A}

\newcommand{\C} {\mathcal C}
\newcommand{\D} {\mathcal D}

\newcommand{\Sc}{\mathcal S}

\newcommand{\End}{\text{\rm End}}
\newcommand{\Ext}{\text{\rm Ext}}
\newcommand{\Hom}{\text{\rm Hom}}

\newcommand{\iso}{\cong}
\newcommand{\Id}{\text{\rm Id}}
\newcommand{\Mod}{\text{\sf Mod-}}

\newcommand{\T}{\mathcal{T}}
\newcommand{\tensor}{\otimes}

\newcommand{\mR}{\mathbb R}
\newcommand{\mS}{\mathbb S}
\newcommand{\mZ}{\mathbb Z}
\newcommand{\spec}{\mathcal{S}p}
\newcommand{\sm}{\wedge}
\renewcommand{\to}{\longrightarrow}

\DeclareMathOperator{\Ho}{Ho}
\DeclareMathOperator{\HoC}{Ho(\C)}

\DeclareMathOperator{\colim}{colim}

\numberwithin{equation}{section}
\newtheorem{theorem}[equation]{Theorem}

\newtheorem{prop}[equation]{Proposition}

\theoremstyle{definition}
\newtheorem{defn}[equation]{Definition}

\newtheorem{rk}[equation]{Remark}
\newtheorem{eg}[equation]{Example}
\newtheorem{egs}[equation]{Examples}

\begin{document}
\title[Morita theory]
{Morita theory in abelian, derived\\ and stable model categories}
\date{\today}
\author{Stefan Schwede}
\address{SFB 478 Geometrische Strukturen in der Mathematik, 
Westf\"alische Wilhelms-Universit\"at M\"unster, Germany}
\email{sschwede@math.uni-muenster.de}

\maketitle

These notes are based on lectures given at the
Workshop on {\em Structured ring spectra and their applications}.
This workshop took place January 21-25, 2002, at the University of Glasgow 
and was organized by Andy Baker and Birgit Richter.

\tableofcontents

\section{Introduction}\label{se:intro}

The paper~\cite{morita} by Kiiti Morita seems to be the 
first systematic study of equivalences between module categories.
Morita treats both contravariant equivalences 
(which he calls {\em dualities} of module categories) 
and covariant equivalences (which he calls {\em isomorphisms} 
of module categories) and shows that they always arise from
suitable bimodules, either via contravariant hom functors
(for `dualities') or via covariant hom functors and tensor products
(for `isomorphisms'). The term `Morita theory' is now used for 
results concerning equivalences of various kinds of module categories.
The authors of the obituary article \cite{Morita obituary} 
consider Morita's theorem 
``probably one of the most frequently used single results in modern algebra''.

In this survey article, we focus on the covariant form of Morita theory,
so our basic question is:

\centerline{When do two `rings' have `equivalent' module categories ?}
\medskip

We discuss this question in different contexts:
\begin{itemize}
\item (Classical) When are the module categories of two rings
equivalent as categories~?
\item (Derived) When are the derived categories of two rings
equivalent as triangulated categories~?
\item (Homotopical) When are the module categories of two ring spectra
Quillen equivalent as model categories~?
\end{itemize}

There is always a related question, which is in a sense more general:

\centerline{What characterizes the category of modules over a `ring' ?}
\medskip

The answer is, mutatis mutandis, always the same:  modules over a `ring'
are characterized by the existence of a `small generator', which
plays the role of the free module of rank one. 
The precise meaning of `small generator' depends on the context,
be it an abelian category, a derived category or a stable model category.
We restrict our attention to categories which have a single small generator; 
this keeps things simple, while showing the main ideas.
Almost everything can be generalized to categories (abelian,
derived or stable model categories) with a {\em set of small generators}.
One would have to talk about {\em ringoids} 
(also called {\em rings with many objects}) and their differential graded
and spectral analogues.

{\bf Background:} for a historical perspective on Morita's
work we suggest a look at the obituary article \cite{Morita obituary}
by Arhangel'skii, Goodearl, and Huisgen-Zimmermann.
The history of Morita theory for derived categories and
`tilting theory' is summarized in Section 3.1 of the book 
by K{\"o}nig and Zimmermann \cite{koenig-zimmermann}. Both sources
contain lots of further references.

For general background material on
derived and triangulated categories, see~\cite{SGA 4 1/2}
(Appendix by Verdier),
\cite{Gelfand-Manin}, \cite{verdier}, or~\cite{Weibel:1994a}.
We freely use the language of model categories, alongside with
the concepts of Quillen adjoint pair and Quillen equivalence.
For general background on model categories 
see Quillen's original article~\cite{Q},
a modern introduction~\cite{DS}, or~\cite{hovey-book} for
a more complete overview.

{\bf Acknowledgments:} The Morita theory in stable model categories
which I describe in Section \ref{sec-Morita stable model categories}
is based on joint work with Brooke Shipley spread over many years 
and several papers; I would like to take this opportunity 
to thank her for the pleasant and fruitful collaboration. 
I would also like to thank Andy Baker and Birgit Richter for organizing the
wonderful workshop {\em Structured ring spectra and their applications}
in Glasgow.

\section{Morita theory in abelian categories}\label{sec-classical Morita}

To start, we review the covariant Morita theory for 
modules; this is essentially the content of Section 3 of~\cite{morita}.
This and related material is treated in more detail 
in~\cite[II \S 3]{bass-Ktheory}, ~\cite[\S 22]{anderson-fuller}
or \cite[\S 18]{Lam}.

\begin{defn}\label{def- abelian small/generator}
Let $\A$ be  an abelian category with infinite sums.
An object $M$ of $\A$ is {\em small} 
if the hom functor $\A(M,-)$ preserves sums; 
$M$ is a {\em generator} if every object of $\A$ 
is an epimorphic image of a sum of 
(possibly infinitely many) copies of $M$.
\end{defn}

The emphasis in the smallness definition is on {\em infinite} sums;
finite sums are isomorphic to finite products, so they are 
automatically preserved by the hom functor. 
For modules over a ring, smallness is closely related to finite generation:
every finitely generated module is small and for projective modules, 
`small' and `finitely generated' are equivalent concepts. 
Rentschler~\cite{rentschler} gives an example of a small module which
is not finitely generated.

A generator can equivalently be defined by the property that the
functor $\A(M,-)$ is faithful, compare~\cite[II Prop.~1.1]{bass-Ktheory}.
A small projective generator is called a {\em progenerator}.
The main example is when  $\A=\Mod R$ is the category of right modules
over a ring $R$. Then the free module of rank one is a small projective
generator.
In this case, a general $R$-module $M$ is a generator for $\Mod R$ if
and only if the free $R$-module of rank one is an epimorphic image
of a sum of copies of $M$.

Here is one formulation of the classical Morita theorem for rings:

\begin{theorem} \label{thm-classical Morita}
For two rings $R$ and $S$, the following conditions are
equivalent. 
\begin{enumerate}
\item[(1)] The categories of right $R$-modules and right $S$-modules are
equivalent.
\item[(2)] The category of right $S$-modules has a small projective generator
whose endomorphism ring is isomorphic to $R$.
\item[(3)] There exists an $R$-$S$-bimodule $M$ such that the functor
\[ - \tensor_R M \ : \ \Mod R \ \to \ \Mod S \] 
is an equivalence of categories.
\end{enumerate}
If these conditions hold, then $R$ and $S$ are said to be
{\em Morita equivalent}.\end{theorem}

Here are some elementary remarks on Morita equivalence.
Condition (1) above is symmetric in $R$ and $S$. So if
an $R$-$S$-bimodule $M$ realizes an equivalence of module categories,
then the inverse equivalence is also realized by an $S$-$R$-bimodule $N$.
Since the equivalences are inverse to each other,
$M \tensor_S N$ is then isomorphic to $R$ as an $R$-bimodule and
$N \tensor_R M$ is isomorphic to $S$ as an $S$-bimodule.
Moreover, $M$ and $N$ are then projective as right modules, and
$N$ is isomorphic to $\Hom_{S}(M,S)$ as a bimodule.

If $R$ is Morita equivalent to $S$, then
the opposite ring $R^{op}$ is Morita equivalent to $S^{op}$.
Indeed, suppose the $R$-$S$-bimodule $M$ and the $S$-$R$-bimodule $N$
satisfy
\[ M \tensor_S N \ \iso \ R \text{\quad and\quad}
N \tensor_R M \ \iso \ S \] 
as bimodules. Since the category of right $R^{op}$-modules is
isomorphic to the category of left $R$-modules,
we can view $M$ as an $S^{op}$-$R^{op}$-bimodule
and $N$ as an $R^{op}$-$S^{op}$-bimodule, and then they provide the
equivalence of categories between $\Mod R^{op}$ and  
$\Mod S^{op}$.
Similarly, if $R$ is Morita equivalent to $S$ and
$R'$ is Morita equivalent to $S'$, then $R\tensor S$ 
is Morita equivalent to $R'\tensor S'$. Here, and in the rest of the paper,
undecoratd tensor products are taken over $\mZ$.

Invariants which are preserved under Morita equivalence include all
concepts which can be defined from the category of modules without reference
to the ring. Examples are the number of isomorphism classes 
of projective modules, of simple modules or of indecomposable modules,
or the algebraic $K$-theory of the ring.
The  center $Z(R)=\{r\in R \, | \, rs=sr \text{ for all } s\in R \}$
is also Morita invariant, since the center of $R$ is isomorphic
to the endomorphism ring of the identity functor of $\Mod R$.
A ring isomorphism 
\[ Z(R) \ \to \ \End(\Id_{\Mod R}) \]
is obtained as follows:
if $r\in Z(R)$ is a central element, then for every $R$-module $M$,
multiplication by $r$ is $R$-linear. So the collection of $R$-homomorphisms
$\{\times r:M\to M\}_{M\in\Mod R}$ is a natural transformation from the
identity functor to itself. 
For more details, see~\cite[II Prop~2.1]{bass-Ktheory}
or \cite[Remark 18.43]{Lam}. In particular, if two {\em commutative} 
rings are Morita equivalent, then they are already isomorphic.

There is a variation of the Morita theorem \ref{thm-classical Morita}
relative to a commutative ring $k$, with essentially the same proof. 
In this version $R$ and $S$ are $k$-algebras, condition (1)
refers to a $k$-linear equivalence of module categories,
condition (2) requires an isomorphism of $k$-algebras 
and in part (3), $M$ has to be a {\em $k$-symmetric} bimodule,
i.e., the scalars from the ground ring $k$ act in the same way from the
left (through $R$) and from the right (through $S$).

We sketch the {\bf proof of Theorem \ref{thm-classical Morita}}
because it serves as the blueprint for analogous results
in the contexts of differential graded rings and ring spectra.

Suppose (1) holds and let
\[ F \ : \ \Mod R \ \to \ \Mod S \] 
be an equivalence of categories. The free $R$-module of rank one is 
a small projective generator of the category of $R$-modules. 
Being projective, small or a generator are categorical conditions, 
so they are preserved by an equivalence of categories.
So the $S$-module $FR$ is a small projective generator 
of the category of $S$-modules. Since $F$ is an equivalence of categories,
it is in particular an additive fully faithful functor.
So $F$ restricts to an isomorphism of rings
\[ F \ : \ R \ \iso \ \End_{R}(R) \ \xrightarrow{\quad \iso\quad} \
 \End_{S}(FR) \ . \]

Now assume condition (2) and let $P$ be a small 
projective $S$-module which generates the category $\Mod S$. 
After choosing an isomorphism $f:R\iso\End_S(P)$, we can view
$P$ as an $R$-$S$-bimodule by setting $r\cdot x=f(r)(x)$
for $r\in R$ and $x\in P$.
We show that $P$ satisfies the conditions of (3) by showing that
the adjoint functors $-\tensor_RP$ and $\Hom_S(P,-)$ are
actually inverse equivalences.

The adjunction unit is the $R$-linear map
\[ X \ \to \   \Hom_{S}(P,X\tensor_RP) \ , \quad
x \longmapsto (y \longmapsto x\tensor y) \ . \]
For $X=R$, the map adjunction unit
coincides with the isomorphism $f$, so it is bijective.
Since $P$ is small, source and target
commute with sums, finite or infinite, so 
the unit is bijective for every free $R$-module.
Since $P$ is projective over $S$, both sides of the adjunction unit
are right exact as functors of $X$. Every $R$-module
is the cokernel of a morphism between free $R$-modules, so
the adjunction unit is bijective in general.

The adjunction counit is the $S$-linear evaluation map
\[ \Hom_{S}(P,Y)\tensor_R P \ \to \ Y \ , \quad
\phi \tensor x \longmapsto \phi(x) \ . \]
For $Y=P$, the counit is an isomorphism since the right
action of $R$ on $\Hom_{S}(P,P)$ arises from the isomorphism
$R\iso\Hom_{S}(P,P)$. Now the argument proceeds 
as for the adjunction unit:
both sides are right exact and preserves sums, finite or infinite, 
in the variable $Y$. Since $P$ is a generator, every $S$-module 
is the cokernel of a morphism between direct sums of copies of $P$, so
the counit is bijective in general.

Condition (1) is a special case of (3), so this finishes 
the proof of the Morita theorem.

\begin{eg}\label{ex-classical matrix}
The easiest example of a Morita equivalence involves matrix algebras.
Any free $R$-module of finite rank $n\geq 1$ is a small projective generator
for the category of right $R$-modules. The endomorphism ring 
\[ \End_{R}(R^n) \ \iso \ M_n(R) \] 
is the ring of $n\times n$ matrices with entries in $R$.
So $R$ and the matrix ring $M_n(R)$ are Morita equivalent.
The bimodules which induce the equivalences of module categories
can both be taken to be $R^n$, but viewed as `row vectors' 
(or $1\times n$ matrices) and `column vectors' (or $n\times 1$ matrices) 
respectively.

Matrix rings do not provide the most general kind of Morita equivalences,
as the example below shows. However, every ring Morita equivalent to $R$
is isomorphic to a ring of the form $eM_n(R) e$ where
$e\in M_n(R)$ is a {\em full idempotent} in the $n\times n$ matrix ring,
i.e., we have $e^2=e$ and $ M_n(R)e  M_n(R)=M_n(R)$.
Indeed, if $P$ is a small projective generator for 
a ring $R$, then $P$ is a summand of a free module of finite rank $n$, say.
Thus  $P$ is isomorphic to the image of an idempotent $n\times n$ matrix
$e$, and then $\End_{R}(P)\iso eM_n(R) e$ as rings.
\end{eg}

\begin{eg} The following example of a Morita equivalence which is not 
of matrix algebra type was pointed out to me 
by M.~K{\"u}nzer and N.~Strickland.
Consider a {\em commutative} ring $R$ and an {\em invertible} $R$-module $Q$.
In other words, there exists another $R$-module $Q'$ and an isomorphism
of $R$-modules $Q\tensor_RQ'\iso R$. Then tensor product with $Q$ over
$R$ is a self-equivalence of the category of right $R$-modules
(with quasi-inverse the tensor product with $Q'$).
This self-equivalence is not isomorphic to the identity functor 
unless $R$ is free of rank one.

Because tensor product with an invertible module $Q$ is an equivalence of
categories, it follows that $Q$ is a progenerator, 
with endomorphism ring isomorphic to $R$.
Moreover, the `inverse' module $Q'$ is isomorphic to the $R$-linear dual
$Q^*=\Hom_R(Q,R)$. Now we consider the direct sum $P=R\oplus Q$, which is
another small projective generator for $\Mod R$.
Then $R$ is Morita equivalent to the endomorphism ring 
of $P$, 
\[ \End_R(P) \ = \ \Hom_R(R\oplus Q, R\oplus Q) \ . \]
As an $R$-module, $\End_R(P)$ is thus isomorphic to 
$R\oplus Q\oplus Q^*\oplus R$. So if $Q$ is not free, then
$\End_R(P)$ is not free over its center, hence not a matrix algebra.

For a specific example we consider the ring 
\[ R \ = \ \mZ[u]/(u^2-5u) \ . \]
We set $Q=(2,u)\lhd R$, the ideal generated by $2$ and $u$.
Then $Q$ is not free as an $R$-module, 
but it is invertible because the evaluation map
\[ \Hom_R(Q,R) \tensor _R Q \ \to \ R \ , 
\quad \phi \tensor x \longmapsto \phi(x) \] 
is an isomorphism.
Note that the inclusion $Q\to R$ becomes an isomorphism after inverting 2;
so after inverting 2 the module $P=R\oplus Q$ is free of rank 2 
and hence the ring $\End_R(P)[\frac{1}{2}]$ is isomorphic to the ring 
of $2\times 2$ matrices over $R[\frac{1}{2}]$.
\end{eg}

The implication (2)$\Longrightarrow$(1) in the Morita theorem
\ref{thm-classical Morita} can be stated in a more general
form, and then it gives a characterization of module categories
as the cocomplete abelian categories with a small projective generator.

\begin{theorem}\label{thm-Gabriel/Michell}
Let $\A$ be an abelian category with infinite sums and a small projective
generator $P$. Then the functor
\[ \A(P,-) \ \colon \ \A \ \to \ \Mod\End_{\A}(P) \]  
is an equivalence of categories.
\end{theorem}
\begin{proof}
We give the same proof as in Bass' book~\cite[II Thm.~1.3]{bass-Ktheory}.
Let us say that an object $X$ of $\A$ is {\em good} if the map
\begin{equation}\label{third eval map}
\A(P,-) \ : \ \A(X,Y) \ \to \ \Hom_{\End_{\A}(P)}(\A(P,X),\A(P,Y)) \end{equation} 
is bijective for every object $Y$ of $\A$. We note that:
\begin{itemize} 
\item The generator $P$ is good since $\A(P,P)$ is the free $\End_{\A}(P)$-module 
of rank one.
\item The class of good objects is closed under sums, finite or infinite:
since $P$ is small, $\A(P,-)$ preserves sums and both sides of
the map \eqref{third eval map} take direct sums in $X$ to direct products.
\item If $f:X\to X'$ is a morphism between good objects in $\A$, 
then the cokernel of $f$ is also good. 
This uses that $P$ is projective and so that $\A(P,-)$ is an exact functor
and both sides of the map \eqref{third eval map} are right exact in $X$.
\end{itemize}
Since $P$ is a generator, every object can be written
as the cokernel of a morphism between sums of copies of $P$.
So every object of $\A$ is good, which precisely means that 
the hom functor $\A(P,-)$ is full and faithful.

It remains to check that every $\End_{\A}(P)$-module is isomorphic 
to a module in the image of the functor $\A(P,-)$. 
The free $\End_{\A}(P)$-module of rank one is the image of $P$.
Since  $\A(P,-)$ commutes with sums, every free module is in the image,
up to isomorphism. Finally, every $\End_{\A}(P)$-module $X$
has a presentation, so it occurs in an exact sequence of $\End_{\A}(P)$-modules
\[ \bigoplus_I \End_{\A}(P) \ \xrightarrow{g} \  \bigoplus_J\End_{\A}(P) \ \to \ X \ \to
\ 0 \ . \]
Since  $\A(P,-)$ is full, the homomorphism $g$ is isomorphic to $\A(P,f)$ 
for some morphism $f:\bigoplus_IP\to\bigoplus_JP$ in $\A$.
Since the functor $\A(P,-)$ is exact, $X$ is the image of the cokernel of $f$.
Thus $\A(P,-)$ is an equivalence of categories. 
\end{proof}

Theorem \ref{thm-Gabriel/Michell} can be applied to the abelian category
of right modules over a ring $S$; then we conclude that for every 
small projective generator $P$ of $\Mod S$ the functor
\[  \Hom_{S}(P,-) \ \colon \ \Mod S \ \to \ \Mod \End_{S}(P) \]  
is an equivalence of categories.
This shows again that condition (2) in the Morita theorem
\ref{thm-classical Morita} implies condition (1).

\section{Morita theory in derived categories}
\label{sec-derived Morita}

Morita theory for derived categories is about the question: 

\centerline{When are the derived categories
$\D(R)$ and $\D(S)$ of two rings $R$ and $S$ equivalent ?}

\medskip

Here the derived category $\D(R)$ is defined from ($\mZ$-graded and unbounded)
chain complexes of right $R$-modules by formally inverting the
quasi-isomorphisms, i.e., the chain maps which induce isomorphisms
of homology groups. Of course, if $R$ and $S$ are Morita equivalent,
then they are also derived equivalent. But it turns out that
derived equivalences happen under more general circumstances.

Rickard \cite{rickard1,rickard2} developed a Morita theory 
for derived categories based on the notion of a tilting complex.
Rickard's theorem did not come out of the blue, and he had built on
previous work of several other people on {\em tilting modules}.
Section 3.1 of the book by K\"onig and Zimmermann \cite{koenig-zimmermann}
gives a summary of the history in this area;
this book also contains many more details, examples and references on the use
of derived categories in representation theory.

We follow Keller's approach from \cite{keller-derivingDG},
based on the (differential graded) endomorphism ring of a tilting complex.
A similar approach to and more applications of
Morita theory in derived categories can be found in the paper~\cite{DG} 
by Dwyer and Greenlees.

\subsection{The derived category}

In this section, $R$ is any ring. All chain complexes are $\mZ$-graded 
and homological, i.e., the differential decreases the degree by 1.

\begin{defn}\label{def-D(R)} A chain complex $C$ of $R$-modules is 
{\em cofibrant} if there exists an exhaustive 
increasing filtration by subcomplexes
\[ 0=C^0 \ \subseteq \ C^1 \ \subseteq \ \cdots \ \subseteq \ C^n 
\ \subseteq \cdots \] 
such that each subquotient $C^n/C^{n-1}$ consists of projective modules
and has trivial differential.
The {\em (unbounded) derived category} $\D(R)$ of the ring $R$ has as objects
the cofibrant complexes of $R$-modules and as morphisms the
chain homotopy classes of chain maps.
\end{defn}

Our definition of the derived category is different from
the usual one. The more traditional way is to start with the homotopy
category of all complexes, not necessarily cofibrant; 
then one uses a calculus of fractions to formally invert the
class of quasi-isomorphisms. These two ways of constructing $\D(R)$ 
lead to equivalent categories.

The {\em shift functor} in $\D(R)$ is given by shifting a complex, i.e.,
\[ (A[1])_n \ = \ A_{n-1} \] 
with differential $d:(A[1])_n=A_{n-1}\to A_{n-2}=(A[1])_{n-1}$
the {\em negative} of the differential of the original complex $A$.
The {\em mapping cone} $C\varphi$ of a chain map $\varphi:A\to B$ is defined by
\begin{equation}\label{eq-def mapping cone}
(C\varphi)_n \ = \ B_n \, \oplus \, A_{n-1} \ , 
\quad d(x,y) \ = \ (dx + \varphi(y),-dy) \ .    \end{equation}
The mapping cone comes with an inclusion $i:B\to C\varphi$ and a projection
$p:C\varphi\to A[1]$ which induce an isomorphism $(C\varphi)/B\iso A[1]$;
if $A$ and $B$ are cofibrant, then so are the shift $A[1]$ 
and the mapping cone.

\begin{rk}\label{rk-cofibrant remark}
The following remarks are meant to give a better feeling for the notion
of `cofibrant complex' and the unbounded derived category.
\begin{itemize}
\item[(i)] The concept of a `cofibrant complex' is closely related to,
but stronger than, a complex of projective modules.  
Indeed, if $C$ is a cofibrant complex, then in every dimension $k\in\mZ$, 
each subquotient $C^n_k/C^{n-1}_k$ is projective.
So $C^n_k$ splits as the sum of the subquotients,
\[ C^n_k \ \iso \ \bigoplus_{i=1}^{n} \ C^i_k/C^{i-1}_k \ . \]
Since $C_k$ is the union of the submodules $C^n_k$, the module
$C_k$ also splits as the sum of the countably many subquotients 
$C^i_k/C^{i-1}_k$. In particular, $C_k$ is a sum of projective modules.
Hence every cofibrant complex is dimensionwise projective.
If $C$ is a complex of projective modules which 
is {\em bounded below}, then it is also cofibrant. 
For example, if $C$ is trivial in negative dimensions, then as the filtration 
we can simply take the (stupid) truncations of $C$, i.e.,
\[ C^i_n \ = \ \begin{cases}
C_n & \text{\ for $n< i$}\\
0   & \text{\ for $n\geq i$.}
\end{cases}\]
So for bounded below complexes, `cofibrant' is equivalent 
to `dimensionwise projective'.

On the other hand, not every complex 
which is dimensionwise projective is also cofibrant. 
The standard example is the complex $C$ 
in which $C_k$ is the free $\mZ/4$-module of rank one
for all $k\in \mZ$, and where every differential $d:C_k\to C_{k-1}$
is multiplication by 2. 
\item[(ii)] Every quasi-isomorphism between cofibrant complexes 
is a chain homotopy equivalence.  
For every complex of $R$-modules $X$, there is a cofibrant
complex $X^{\text{c}}$ and a quasi-isomorphism $X^{\text{c}}\to X$;
together these two facts essentially prove that the
derived category $\D(R)$ enjoys the universal property of the localization
of the category of chain complexes of $R$-modules with the class 
of quasi-isomorphisms inverted.
These properties are very analogous to the properties that CW-complexes
have among all topological spaces:
every weak equivalences between CW-complexes is a homotopy equivalence
and every space admits a CW-approximation.
This analogy is made precise in \cite[Part III]{kriz-may}.
\item[(iii)]  The concept of a cofibrant chain complex is closely related to 
that of a {\em K-projective} complex
as defined by Spaltenstein \cite[Sec.\ 1.1]{spaltenstein} (who attributes
this notion to J.~Bernstein; Keller \cite[8.1.1]{koenig-zimmermann}
calls this {\em homotopically projective}).
A chain complex is K-projective if every chain map into 
an acyclic complex (i.e., a complex with trivial homology) 
is chain null-homotopic. 

Every cofibrant complex is K-projective. Conversely, every K-projective
complex $X$ is chain homotopy equivalent to a cofibrant complex.
Indeed, we can choose a cofibrant replacement, i.e., a cofibrant
complex $X^c$ and a quasi-isomorphism $q:X^c\to X$; 
the mapping cone $Cq$ is then acyclic. We have a short exact sequence
of chain homotopy classes of chain maps in $Cq$,
\[ [X^c[1],Cq] \ \to \ [Cq,Cq] \ \to \ [X,Cq] \ ; \]
both $X$ and $X^c$ are K-projective, so the left and
right groups are trivial. Thus the identity map of $Cq$ is null-homotopic
and so the mapping cone of $q$ is chain contractible. Thus the map
$q:X^c\to X$ is a chain homotopy equivalence.
\item[(iv)]
We use the term `cofibrant' complex because they are the cofibrant objects
in the {\em projective} model category structure on chain complexes
of $R$-modules \cite[2.3.11]{hovey-book}. In this model structure,
the weak equivalences are the quasi-isomorphisms, the fibrations are the
surjections and the cofibrations are the injections whose cokernel is
cofibrant in the sense of Definition \ref{def-D(R)}.
In particular, every chain complex is fibrant 
in the projective model structure.

Since every object is fibrant and because the fibrations (= surjections)
and weak equivalences (= quasi-iso\-mor\-phisms) already 
have well-established names, it is unnecessary, and overly complicated,
to use the language of model categories in order to work 
with the derived category of a ring.
\item[(v)] There is a `dual' approach to the derived category $\D(R)$ 
as the homotopy category of `fibrant' complexes
(attention: these are {\em not} the fibrant objects in the
projective model structure -- there every complex is fibrant).
This uses the notion of a `K-injective' \cite[Sec.\ 1.1]{spaltenstein}
or `homotopically injective' \cite[8.1.1]{koenig-zimmermann}
complex, or the {\em injective model structure}
on the category of complexes of $R$-modules~\cite[2.3.13]{hovey-book}. 
It is often useful to have both descriptions available. 
Given arbitrary chain complexes $C$ and $D$,
we choose a cofibrant/K-projective resolution 
$C^{\text{c}}\stackrel{\sim}{\to} C$
and a fibrant/K-injective resolution $D\stackrel{\sim}{\to} D^{\text{f}}$. 
Then the maps induce isomorphisms of chain homotopy classes of chain maps
\[ [C^{\text{c}},D] \ \xrightarrow{\ \iso\ } \ [C^{\text{c}},D^{\text{f}}] \ \xleftarrow{\ \iso\ } \
 \ [C,D^{\text{f}}] \ . \] 
\end{itemize}
\end{rk}

There is an additive functor
\[ [0] \ : \ \Mod R \ \to \ \D(R) \] 
which is a fully faithful embedding onto the full subcategory of the
derived category consisting of the complexes whose homology is concentrated 
in dimension zero. So we can think of the $R$-modules as sitting inside
the derived category $\D(R)$. Had we defined the derived category from the
category of all complexes of $R$-modules by formally inverting the
quasi-isomorphisms, then we could define the complex $M[0]$
by putting the $R$-module $M$ in dimension 0, and taking trivial
chain modules everywhere else. With our present definition of $\D(R)$
we let $M[0]$ be a choice of resolution $P_{\bullet}$ 
of $M$ by projective $R$-modules.
Such a resolution is unique up to chain homotopy equivalence 
and it is cofibrant when viewed as a chain complex
(by Remark \ref{rk-cofibrant remark} (i) above), 
Moreover, every $R$-linear map $M\to N$ 
is covered by a unique chain homotopy class between the chosen
projective resolutions. In other words, we really get a functor from
$R$-modules to the derived category $\D(R)$, together with a natural
isomorphism $H_0(M[0])\iso M$.

The usual definition of Ext-groups involves a choice 
of projective resolution  $P_{\bullet}$ of the source module $M$,
and then $\Ext^n_R(M,N)$ can be defined as the chain homotopy classes
of chain maps from the resolution $P_{\bullet}$ to $N$, 
shifted up into dimension $n$.
We get the same result if $N$ is also replaced by a projective resolution;
this says that Ext-groups can be obtained from the derived category via
\begin{equation}\label{Ext in D(R)}
\Ext^n_R(M,N)  \iso \  \D(R)(M[0],N[n]) \ . \end{equation}

The derived category of a ring has more structure. 
The category $\D(R)$ is additive 
since the homotopy relation for chain maps is additive. 
But $\D(R)$ is no longer an abelian category such as the category 
of chain complexes.
Indeed, notions such as `monomorphism', `epimorphisms, `kernels' 
for chain maps do not interact well with the passage to chain homotopy classes.
The {\em distinguished triangles} in $\D(R)$ are what is left of the
abelian structure on the category of chain complexes, and $\D(R)$
is an example of a {\em triangulated category}.

The distinguished triangles are the diagrams which are isomorphic 
in $\D(R)$ to a mapping cone triangle.
More precisely, a diagram in $\D(R)$ of the form
\begin{equation}\label{generic triangle} 
X \ \xrightarrow{f} \ Y \ \xrightarrow{g} \ Z \ \xrightarrow{h} \ X[1] 
\end{equation}
is called a {\em distinguished triangle} if and only if there exists
a chain map  $\varphi:A\to B$ between cofibrant complexes and isomorphisms
$\iota_1:A\iso X$, $\iota_2:B\iso Y$ and $\iota_3:C\varphi\iso  Z$  
in $\D(R)$ such that the diagram 
\[\xymatrix{  A \ar[r]^-{\varphi} \ar[d]_{\iota_1} & 
B \ar[r]^-{i} \ar[d]_{\iota_2} &
C\varphi \ar[r]^{p} \ar[d]^{\iota_3} & A[1] \ar[d]^{\iota_1[1]}\\
X \ar[r]_-{f} & Y \ar[r]_-{g} & Z \ar[r]_-{h} & X[1] }\]
commutes in $\D(R)$.

We do not want to reproduce the complete definition 
of a triangulated category here;
the data of a triangulated category consists of 
\begin{enumerate}
\item[(i)]  an additive category $\T$,
\item[(ii)] a self-equivalence $[1]:\T\to \T$ called the {\em shift} 
functor and
\item[(iii)] a class of {\em distinguished triangles}, i.e., a collection
of diagrams in $\T$ of the form \eqref{generic triangle}.
\end{enumerate}
This data is subject to several axioms which can be found for example in
\cite{verdier}, \cite[Sec.\ 10.2]{Weibel:1994a} \cite{neeman-triangle book} 
or \cite[2.3]{koenig-zimmermann}.

Distinguished triangles are the source 
of many long exact sequences that come up in nature.
Indeed, the axioms which we have suppressed imply in particular that for
every distinguished triangle of the form \eqref{generic triangle} 
and every object $W$ of $\T$ the sequence of abelian morphism groups
\[ \T(W,X) \ \xrightarrow{f_*} \ \T(W,Y) \ \xrightarrow{g_*} \ 
\T(W,Z) \ \xrightarrow{h_*} \ \T(W,X[1]) \]
is exact. One of the axioms also says that one can `rotate' triangles,
i.e., a sequence \eqref{generic triangle} is a distinguished triangle 
if and only if the sequence
\[ Y \ \xrightarrow{\ g\ } \ Z \ \xrightarrow{\ h\ } \ X[1] \ 
\xrightarrow{-f[1]} \ 
Y[1] \] 
is a distinguished triangle. So if we keep rotating a distinguished triangle
in both directions and take morphisms from a fixed object $W$,
we end up with a long exact sequence of abelian groups
\begin{align}\label{generic LES}
\cdots \ \T(W,X)  \xrightarrow{\ f_*\ }  &\T(W,Y)  
\xrightarrow{\ g_*\ } \T(W,Z) \xrightarrow{\ h_*\ } \T(W,X[1]) \\
& \xrightarrow{-f[1]_*} \T(W,Y[1])  \xrightarrow{-g[1]_*}  
\T(W,Z[1])  \xrightarrow{-h[1]_*}  \T(W,X[2]) \ \cdots\nonumber \end{align}
The axioms of a triangulated category also guarantee a similar 
long exact sequence when taking morphism {\em from} a triangle 
(and its rotations) {\em into} a fixed object $W$.

In the derived category $\D(R)$, the long exact sequence \eqref{generic LES}
becomes something more familiar when we take $W=R[0]$, the free $R$-module 
of rank one, viewed as a complex concentrated in dimension zero. 
For every chain complex $C$, cofibrant or not, the chain homotopy
classes of morphisms from  $R$ to $C$ are naturally isomorphic to
the homology module $H_0C$; so the long exact sequence \eqref{generic LES}
specializes to the long exact sequence of homology modules
\[  \cdots \ \to \ H_0(A) \ \xrightarrow{\ H_0(\varphi)\ } 
\ H_0(B) \ \to \ H_0(C\varphi)
\  \xrightarrow{\ \delta \ } \ H_{-1}(A) \ \to \ \cdots  \ . \]

{\bf Small generators.}
We will often require infinite direct sums in a triangulated category. 
The unbounded derived category $\D(R)$ has direct sums, finite and infinite. 
Indeed, the direct sum 
of any number of cofibrant complexes is again cofibrant 
(take the direct sum of the filtrations which are required 
in Definition \ref{def-D(R)}), and this also represents the direct sum
in $\D(R)$. This is one point where it is important to allow unbounded
complexes. There are variants of the derived category 
which start with complexes which are bounded or bounded below.
One also gets triangulated categories in much the same way as for $\D(R)$, 
but for example the countable family $\{R[-n]\}_{n\geq 0}$ has no direct sum
in the bounded or bounded below derived categories.

In the Morita equivalence questions, a suitably defined notion of
`small generator' pops up regularly. 
The following concepts for triangulated categories
are analogous to the ones for abelian categories 
in Definition \ref{def- abelian small/generator}.

\begin{defn}\label{def- triangulated small/generator}
Let $\T$ be a triangulated category with infinite coproducts. 
An object $M$ of $\T$ is {\em small} 
if the hom functor $\T(M,-)$ preserves sums; 
$M$ is a {\em generator} if there is no proper
full triangulated subcategory of $\T$ (with shift  
and triangles induced from $\T$) which contains $M$ and
is closed under infinite sums.
\end{defn}

As in abelian categories, the hom functor $\T(M,-)$ automatically
preserves finite sums. What we call `small' is sometimes called
{\em compact} or {\em finite}
in the literature on triangulated categories.
A triangulated category with infinite coproducts and a set 
of small generators is often called {\em compactly generated}.

The class of small objects in any triangulated category is closed under
shifting in either direction, taking finite sums
and taking direct summands. Moreover, if two of the three objects 
in a distinguished triangle are small, then so is the third one 
(one has to exploit that the morphisms 
from a distinguished triangle into a fixed object give rise 
to a long exact sequence).

There is a convenient criterion for when a {\em small} object $M$
generates a triangulated category $\T$ with infinite coproducts:
$M$ generates $\T$ in the sense of 
Definition \ref{def- triangulated small/generator}
if and only  if it `detects objects',
i.e., an object $X$ of $\T$ is trivial if and only if there are  
no graded maps from $M$ to $X$, i.e.
$\T(M[n], X)=0$ for all $n\in \mZ$.
For the equivalence of the two conditions, see for example
\cite[Lemma 2.2.1]{ss-classification}.

The complex $R[0]$ consisting of the free module of rank one
concentrated in dimension 0 is a small generator 
for the derived category $\D(R)$. 
Indeed, morphisms in $\D(R)$ out of the complex $R[0]$ represent homology, 
i.e., there is a natural isomorphism
\[ \D(R)(R[n],C) \ \iso \ H_n(C) \]
for every cofibrant complex $C$.
Since homology commutes with infinite sums,  the complex $R[0]$
is small in  $\D(R)$.
Moreover, if all the morphism groups $\D(R)(R[n],C)$ are trivial 
as $n$ ranges over the integers, 
then the complex $C$ is acyclic, hence contractible, and
so it is trivial in the derived category $\D(R)$. 
In other words, mapping out of shifted copies of the complex $R[0]$
detects whether an object in $\D(R)$ is trivial or not,
so $R[0]$ is also a generator, by the previous criterion.

There is a nice characterization of the small objects in the derived category
of a ring. Every bounded complex of finitely generated projective modules 
is built from summands of the small object $R[0]$ by shifts 
and extensions in triangles.
Since the class of small objects is closed under these operations,
a bounded complex of finitely generated projective modules is small.
Conversely, these are the only small objects, up
to isomorphism in $\D(R)$:

\begin{theorem}\label{thm-small in D(R)}
Let $R$ be a ring. A complex of $R$-modules is small
in the derived category $\D(R)$ if and only if it is quasi-isomorphic 
to a bounded complex of finitely generated projective $R$-modules.
\end{theorem}

The proof that every small object in $\D(R)$ is quasi-isomorphic 
to a bounded complex of finitely generated projective modules
is more involved.
It is a special case of a result about triangulated categories $\T$ 
with a set of small generators.
Neeman \cite{neeman-TTYBR} showed that every
small object in $\T$ is a direct summand of an iterated extension
of finitely many shifted generators.
The proof can also be found in \cite[5.3]{keller-derivingDG}.

There are non-trivial triangulated categories in which only the  
zero objects are small, see for example~\cite{keller-smashing}
or~\cite[Cor. B.13]{hovey-strickland-kn}.
If a triangulated category has a set of generators, then the  
coproduct of all of them is a single generator.
However, an infinite coproduct of non-trivial small objects is not small.
So the property of having a single small generator is something special.
In fact we see in Theorem \ref{thm-main-one} below that this  
condition characterizes the module
categories over ring spectra among the stable model categories.
A triangulated category need not have a {\em set} of generators
whatsoever (one could consider all objects, but in general these form 
a proper class), for example $K(\mZ)$, the homotopy category
of chain complexes of abelian groups, 
is not generated by a set \cite[E.3]{neeman-triangle book}.

{\bf Equivalences of triangulated categories.}
A functor between triangulated categories is called exact
if it commutes with shift and preserves distinguished triangles.
More precisely, $F\colon \Sc \to \T$ is {\em exact} it is is equipped with
a natural isomorphism 
\mbox{$\iota_X:F(X[1])\iso F(X)[1]$} such  that 
for every distinguished triangle \eqref{generic triangle}
the sequence
\[ F(X) \ \xrightarrow{F(f)} \ F(Y) \ \xrightarrow{F(g)} \ F(Z) 
\ \xrightarrow{\iota_X\circ F(h)} \ F(X)[1]  \]
is again a distinguished triangle.
An exact functor is automatically additive.
An {\em equivalence of triangulated categories} is an equivalence of categories
which is exact and whose inverse functor is also exact.

Exact equivalences between derived categories preserve
all concepts which can be defined from $\D(R)$ using only
the triangulated structure.
One such invariant is the Grothendieck group $K_0(R)$,
defined as the free abelian group generated by the isomorphism classes
of finitely generated projective $R$-modules, modulo the relation
\[[P] \ + \ [Q]  \ = \  [P\oplus Q] \ . \]

For any compactly generated triangulated category $\T$, 
the Grothendieck group $K_0(\T)$ is defined 
as the free abelian group generated by the isomorphism classes
of small objects in $\T$, modulo the relation
\[ [X] \ + \ [Z] \ = \ [Y] \] 
for every distinguished triangle
\[ X \ \to \ Y \ \to \ Z \ \to \ X[1] \]
involving small objects $X$, $Y$ and $Z$
(the morphisms in the triangle do not affect the relation).
The split triangle
\[   X \ \xrightarrow{(1,0)} \ X\oplus Y \ 
\xrightarrow{0\choose 1} \ Y \ \xrightarrow{\ 0 \ } \ X[1] \] 
is always distinguished, so the relation $[X\oplus Y]=[X]+ [Y]$ 
holds in $K_0(\T)$. So for every ring $R$, the assignment
\begin{equation}\label{K(R) to K(D(R))} 
K_0(R) \ \to\ K_0(\D(R)) \ , \ \quad [P] \ \longmapsto \ [P[0]] \end{equation}
defines a group homomorphism. This is in fact an isomorphism,
see~\cite[Sec.\ 7]{grothendieck-groupes des classes}. 
The inverse takes the class  in $K_0(\D(R))$ of a bounded complex $C$ of 
finitely generated projective modules to its `Euler characteristic',
\[ \sum_{n\in\mZ} \, (-1)^n[C_n] \ \in \ K_0(R) \ . \]

It is much less obvious that constructions such as the
center of a ring, Hochschild and cyclic homology and
the higher Quillen $K$-groups are also invariants of the derived category. 
In contrast to the Grothendieck group $K_0$,
there is no construction which produces these groups from the
triangulated structure of $\D(R)$ only. The proof that two derived equivalent
rings share these invariants uses 
the `tilting theory' which we outline in the next section.
More precisely, if  $R$ and $S$ are derived equivalent 
flat algebras over some commutative ground ring,
then there exists a {\em two-sided tilting complex},
i.e., a chain complex $C$ of $R$-$S$-bimodules such that the
functor $-\tensor_RC$ induces a (possibly different) 
derived equivalence~\cite{rickard2}. 
Tensor product with the bimodule complex $C$ then induces an 
equivalence of $K$-theory spaces by the work of 
Thomason-Trobaugh \cite[Thm. 1.9.8]{Thomason-Trobaugh}.
Without the flatness assumption, the Waldhausen categories  
of small, cofibrant chain complexes can be related 
through an intermediate category of differential graded modules,
to still obtain an equivalence of $K$-theory spaces;
for more details we refer to~\cite{dugger-shipley}. 
Similarly, Hochschild homology and cohomology
(see~\cite[Prop. 2.5]{rickard2}, or, including 
the Gerstenhaber bracket, see~\cite{keller-HochschildPicard})
and cyclic homology (see~\cite{keller-cyclic invariance,
keller-cyclic invariance and localization}) are isomorphic for
derived equivalent rings which are flat algebras over some 
commutative ground ring.
The invariance of the center under derived equivalence
is established in~\cite[Prop.\ 9.2]{rickard1} 
or~\cite[Prop. 6.3.2]{koenig-zimmermann}.

\medskip

There is a certain general argument which we will use several times 
to verify that certain triangulated functors are equivalences,
so we state it as a separate proposition.
This Proposition \ref{triangulated comparison prop} is a version 
of `Beilinson's Lemma'~\cite{beilinson} and is typically applied
when $F$ is the total derived functor of a suitable left adjoint.
In the following proposition, it is crucial that the functor
$F$ be defined and exact on the entire triangulated category $\Sc$.
It is easy to find non-equivalent triangulated categories $\Sc$ and $\T$ 
with infinite sums and small generators $P$ and $Q$ respectively such that
\[  \Sc(P,P)_* \ \iso \ \T(Q,Q)_* \]
as graded rings. For example one can take a differential graded ring $A$
with a non-trivial triple Massey product and consider
derived categories $\Sc=\D(A)$ and $\T=\D(H^*A)$ (where the
cohomology ring of $A$ is given the trivial differential).

\begin{prop}\label{triangulated comparison prop}
Let $F\colon \Sc \to \T$ be an exact functor
between triangulated categories with infinite sums. Suppose that
$F$ preserves infinite sums and $\Sc$ has a small generator $P$ such that 
\begin{enumerate}
\item[(i)] $FP$ is a small generator of $\T$ and
\item[(ii)] for all integers $n$, the map
\[ F\ \colon \ \Sc(P[n],P) \ \to \ \T(FP[n],FP) \] 
is bijective
\end{enumerate}
Then $F$ is an equivalence of categories.
\end{prop}
\begin{proof}
We consider the full subcategory of $\Sc$ consisting of those $Y$ 
for which the map
\begin{equation}\label{eq-EQ1} 
F \ : \  \Sc(P[n],Y) \ \to \ \T(FP[n],FY) \end{equation}
is bijective for all $n\in\mZ$. By assumption this subcategory contains $P$.
Since $F$ is exact, the subcategory is closed under extensions. 
Since $P$ and $FP$ are small and $F$ preserves coproducts,
this subcategory is also closed under coproducts.
Since $P$ generates $\Sc$, the map \eqref{eq-EQ1} 
is thus bijective for arbitrary $Y$.

Similarly for arbitrary but fixed $Y$
the full subcategory of $\Sc$ consisting of those $X$ for which the map
$F:\Sc(X,Y)\to \T(FX,FY)$ is bijective is closed under extensions
and coproducts. By the first part, it also contains $P$, 
so this subcategory is all of $\Sc$.
In other words, $F$ is full and faithful.

Now we consider the full subcategory of $\T$ of objects which are
isomorphic to an object in the image of $F$.
This subcategory contains the generator $FP$ and
it is closed under shifts and coproducts since these are preserved by $F$.
We claim that this subcategory is also closed under extensions.
Since $FP$ generates $\T$, this shows that $F$ is essentially surjective
and hence an equivalence.

To prove the last claim we consider a distinguished triangle
\[ X \ \xrightarrow{\ f\ } \ Y \ \to \ Z \ \to \ X[1] \ . \]
Since the subcategory under consideration is closed under isomorphism 
and shift in either direction we can assume that $X=F(X')$ and $Y=F(Y')$ 
are objects in the image of $F$.
Since $F$ is full there exists a map $f':X'\to Y'$ satisfying $F(f')=f$. 
We can then choose a  mapping cone for the map $f'$ and a compatible map 
from $Z$ to $F(\mbox{Cone}(f'))$ which is necessarily an isomorphism. 
\end{proof}

\subsection{Derived equivalences after Rickard and Keller}
\label{sec-Rickard's tilting theorem}

In this section we state and prove Rickard's 
``Morita theory for derived categories''.
Rickard shows in~\cite[Thm.\ 6.4]{rickard1} 
that the existence of a {\em tilting complex}
is necessary and sufficient for an equivalence between the
unbounded derived categories of two rings.
A tilting complex is a special small generator of the derived category,
see Definition \eqref{def-tilting complex} below.
The idea to use differential graded algebras in the proof 
is due to Keller~\cite{keller-derivingDG},  
and we closely follow his approach.

The notion of a tilting complex comes up naturally when we examine
the properties of the preferred generator $R[0]$
of the derived category $\D(R)$.
First of all, the free $R$-module of rank one, considered as a complex
concentrated in dimension zero, is a small generator 
of the derived category $\D(R)$.
Since $R$ is a free module, it has no self-extensions. 
Because Ext groups can be identified with morphisms in the derived category
(see \eqref{Ext in D(R)}), this means that the graded self-maps of the
complex $R[0]$ are concentrated in dimension zero:
\[ \D(R)(R[n],R) \ = \ 0 \text{ \ for $n\neq 0$.} \]
A tilting complex is any complex which also has these properties. 
Hence the definition is made so that the image of $R[0]$ 
under an equivalence of triangulated categories is a tilting complex.

\begin{defn}\label{def-tilting complex} 
A {\em tilting complex} for a ring $R$ is a bounded complex $T$ 
of finitely generated projective $R$-modules which generates
the derived category  $\D(R)$ and whose graded ring of self maps
${\D(R)}(T,T)_*$ is concentrated in dimension zero.
\end{defn}

Special kinds of tilting complexes are the {\em tilting modules};
we give examples of tilting modules and tilting complexes in Section 
\ref{sec-tilting examples}.
The following theorem is due to Rickard~\cite[Thm.\ 6.4]{rickard1}.

\begin{theorem} \label{thm-tilting} 
For two rings $R$ and $S$ the following conditions are equivalent.\\
\hspace*{0.5cm}{\em (1)} The unbounded derived categories of $R$ and $S$
are equivalent as triangulated categories.\\
\hspace*{0.5cm}{\em (2)} There is a tilting complex 
$T$ in $\D(S)$ whose endomorphism ring $\D(S)(T,T)$ is isomorphic to $R$. 

Moreover, conditions {\em (1)} and {\em (2)} are implied by the condition\\
\hspace*{0.5cm}{\em (3)} There exists a chain complex 
of $R$-$S$-bimodules $M$ such that the derived tensor product functor 
\[ - \tensor^L_R M \ : \ \D(R) \ \to \ \D(S) \]
is an equivalence of categories.\\

If $R$ or $S$ is flat as an abelian group, then all three conditions
are equivalent.
\end{theorem}

Instead of using the unbounded derived category, one can replace 
condition {(1)} by an equivalence between the full subcategories of 
homologically bounded below or small objects inside the derived categories,
see for example \cite[Thm.\ 6.4]{rickard1}.
There is a version relative to a commutative ring $k$. 
Then $R$ and $S$ are $k$-algebras, conditions (1) and (3) 
then refer to $k$-linear equivalences of derived categories, 
condition (2) requires an isomorphism of $k$-algebras 
and in the addendum, one of $R$ or $S$ has to be flat as a $k$-module.

\begin{rk}
A derived equivalence  $F$ from $\D(R)$ to $D(S)$ 
which is not already a Morita equivalence maps the $R$-modules inside $\D(R)$ 
(i.e., complexes with homology concentrated in dimension 0) 
``transversely'' to the $S$-modules inside $\D(S)$; 
more precisely, for an $R$-module $M$, the
complex $F(M[0])$ can have non-trivial homology in several, or even 
in infinitely many dimensions; we give an example 
in \ref{ex-derived spread out matrix} below.
However, the homology of $F(M[0])$ is always bounded below.

A related point is that we can {\em not} recover the module category $\Mod R$ 
from $\D(R)$, viewed as an abstract triangulated category. 
This is because we cannot make sense of
``complexes with homology concentrated in dimension 0'' unless we specify
a homology functor like $H_0$, or we single out the class of complexes
with homology in non-negative dimensions.
This sort of extra structure is called a {\em t-structures}~\cite[1.3]{BBD} 
on a triangulated category.
Every t-structure has a {\em heart}, an abelian category which plays the
role of complexes in $\D(R)$ with homology concentrated in dimension zero.
\end{rk}

The most involved part of the tilting theorem is the implication 
(2)$\Longrightarrow$(1), i.e., showing that a tilting complex gives rise to a 
derived equivalence. The proof we give is due to Keller; 
in the original paper~\cite{keller-derivingDG}, 
his setup is more general (he works in differential graded categories
in order to allow `many generator' versions). 
In the special case of interest for us,
the exposition simplifies somewhat~\cite[Ch.\ 8]{koenig-zimmermann}.
Given a tilting complex $T$ in $\D(S)$, the comparison
between  the derived categories of $R$ and $S$ passes through 
the derived category of a certain {\em differential graded ring}  
(generalizing the derived category of a an ordinary ring),
namely the
{\em endomorphism DG ring} $\End_S(T)$ of the tilting complex $T$
(generalizing the endomorphism ring).
So we start by introducing these new characters.

\begin{defn}\label{def-DG rings and DG modules}
A {\em differential graded ring} is a $\mZ$-graded ring $A$ together
with a differential $d$ of degree $-1$ which satisfies the Leibniz rule
\begin{equation}\label{Leibniz} 
d(a\cdot b) \ = \ d(a)\cdot b \, + \, (-1)^{|a|} \, a\cdot d(b) \end{equation}
for all homogeneous elements $a,b\in A$.
A {\em differential graded right module} (or {\em DG module} for short)
over a differential graded ring $A$ consists of a graded right $A$-module 
together with a differential $d$ of degree $-1$ 
which satisfies the Leibniz rule \eqref{Leibniz}, 
but where now $a$ is a homogeneous element of the module
and $b$ is a homogeneous element of $A$.
A {\em homomorphism} of DG modules is a homomorphism of graded $A$-modules
which is also a chain map. A {\em chain homotopy} between homomorphisms 
of DG modules is a homomorphism of graded $A$-modules of degree 1
which is also a chain homotopy.

A differential graded $A$-module $M$ is {\em cofibrant}
if there exists an exhaustive increasing filtration by sub DG modules
\[ 0=M^0 \ \subseteq \ M^1 \ \subseteq \ \cdots \ \subseteq \ M^n 
\ \subseteq \cdots \] 
such that each subquotient $M^n/M^{n-1}$ is a direct summand of
a direct sum of shifted copies of $A$.
The {\em derived category} $\D(A)$ of the differential graded ring $A$ 
has as objects the cofibrant DG modules over $A$ and as morphisms the
chain homotopy classes of DG module homomorphisms.
\end{defn}

Up to chain homotopy equivalence, the cofibrant DG modules are the ones
which have Keller's `property (P)' in \cite[3.1]{keller-derivingDG}.
A cofibrant differential graded module is sometimes called `semi-free' 
or a `cell module'~\cite[Part III]{kriz-may} (up to direct summands).

\begin{rk}\label{rk-DG module remarks}
We need some facts about differential graded rings and
modules which are not very difficult to prove, but which we do not want to
discuss in detail.
\begin{itemize}
\item[(i)] Several of the remarks from \ref{rk-cofibrant remark} 
carry over from rings to DG rings.
A cofibrant DG $A$-module is projective as a graded $A$-module, ignoring
the differential.  
Every quasi-isomorphism between cofibrant DG modules
is a chain homotopy equivalence, and every DG module can be approximated
up to quasi-isomorphism by a cofibrant one.
A DG module is called
{\em homotopically projective} if every homomorphism into an acyclic 
DG module is null-homotopic. Then a DG module is homotopically projective 
if and only if it is chain homotopy equivalent, as a DG module, to a cofibrant
DG module.
\item[(ii)] The derived category $\D(A)$ of a differential graded ring $A$
is naturally a triangulated category. The shift functor is again given
by reindexing a DG module, and distinguished triangles arise from 
mapping cones as for the derived category of a ring.
The only thing to note is that for a homomorphism $f:M\to N$ of
DG modules over $A$, the mapping cone becomes a graded $A$-module
as the direct sum $N\oplus M[1]$, and this $A$-action satisfies the
Leibniz rule with respect to 
the mapping cone differential \eqref{eq-def mapping cone}.
\item[(iii)] Suppose that $f:A\to B$ is a homomorphism of differential graded
rings, i.e., $f$ is a multiplicative chain homomorphism.
Then extension of scalars $M\mapsto M\tensor_AB$ 
is exact on cofibrant differential graded modules 
(since the underlying {\em graded} modules over the graded ring underlying $A$
are projective), it takes cofibrant modules to cofibrant modules,
and it preserves the chain homotopy relation.
So extension of scalars induces an exact functor
on the level of derived categories
\begin{equation}\label{derived f_* and f^*} 
\xymatrix{ \D(A)\ \ar@<.3ex>^{\ -\tensor^L_AB\ }[r] & 
\ar@<.3ex>^{\ f^*\ }[l]\  \D(B)} \ ,  \end{equation}
called the {\em left derived functor}.
This derived functor has an exact right adjoint $f^*$
induced by restriction of scalars along $f$.
This is not completely obvious with our definition of the derived category, 
since a cofibrant differential graded $B$-module is usually {\em not} 
cofibrant when viewed as a DG module over $A$ via $f$.

If $f:A\to B$ is a {\em quasi-isomorphism} of differential graded rings,
then the derived functors of restriction and extension of scalars
\eqref{derived f_* and f^*} are inverse {\em equivalences} of
triangulated categories.
\item[(iv)] Suppose $A$ is a differential graded ring whose homology is
concentrated in dimension zero. Then $A$ is quasi-isomorphic,
as a differential graded ring, to the zeroth homology ring $H_0A$.
Indeed, a chain of two quasi-isomorphisms is given by
\[ A \ \xleftarrow{\text{\ inclusion\ }} \ A_+ \ 
\xrightarrow{\text{\ projection\ }} \ H_0(A) \ . \] 
Here  $A_+$ is the differential graded  sub-ring of $A$ given by
\[ (A_+)_n \ = \ \begin{cases}
\qquad \quad A_n & \text{ for $n> 0$} \\
\text{Ker} \left( d:A_0\to A_{-1}\right) & \text{ for $n=0$, and }\\
\qquad \quad 0 & \text{ for $n<0$.} \end{cases} \] 
Since the homology of $A$ is trivial in negative dimensions, the
inclusion $A_+\to A$ is a quasi-isomorphism.
Since $A_+$ is trivial in negative dimensions,
the projection $A_+\to H_0(A_+)$ is a homomorphism of differential graded
rings, where the target is concentrated in dimension zero. 
This projection is also a quasi-isomorphism since the homology of $A$,
and hence that of $A_+$, is trivial in positive dimensions.
\end{itemize}
\end{rk}

{\bf Homomorphism complexes.} 
Let $A$ be a DG ring and let $M$ and $N$ be DG modules over $A$, 
not necessarily cofibrant.
We defined the {\em homomorphism complex} $\Hom_A(M,N)$ as follows. 
In dimension $n\in\mZ$, the chain group $\Hom_A(M,N)_n$ 
is the group of graded $A$-module homomorphisms of degree $n$, i.e.,
\[ \Hom_A(M,N)_n \ = \ \Hom_A(M[n],N) \ . \]
The differentials of $M$ and $N$
do not play any role in the definition of the chain groups, 
but they enter in the formula for the differential 
which makes $\Hom_A(M,N)$ into a chain complex.
This differential $d\colon\Hom_A(M,N)_n\to \Hom_A(M,N)_{n-1}$
is defined by
\begin{equation}\label{eq:d in Hom complex} 
d(f) \ = \  d_N \circ f - (-1)^n f\circ d_M \ . 
\end{equation}
Here $f$ is a graded $A$-module map of degree $n$ and the composites
$d_N \circ f$ and $f\circ d_M$ are then graded $A$-module maps of degree $n-1$.

With this definition, the 0-cycles in $\Hom_A(M,N)$ 
are those graded $A$-module maps $f$ which satisfy
$d_N \circ f - f\circ d_M = 0$,
so they are precisely the DG homomorphisms from $M$ to $N$. 
Moreover, if $f,g:M\to N$ are two graded $A$-module maps,
then the difference $f-g$ is a coboundary in the complex 
$\Hom_A(M,N)$ if and only if $f$ is chain homotopic to $g$. 
So we have established a natural isomorphism
\[ H_0 \left(\Hom_A(M,N)\right) \ \iso  \ [M,N] \] 
between the zeroth homology of the complex $\Hom_A(M,N)$ 
and the chain homotopy classes of DG $A$-homomorphisms from $M$ to $N$. 

Now suppose that we have a third DG module $L$.
Then the composition of graded $A$-module maps gives a bilinear pairing
between the homomorphism complexes 
\[ \circ \ : \ \Hom_A(N,L)_m \ \times \ \Hom_A(M,N)_n
\ \to \ \Hom_A(M,L)_{m+n} \ . \]
Moreover, composition and the differential \eqref{eq:d in Hom complex}
satisfy the Leibniz rule, i.e., for graded $A$-module maps
$f:M\to N$ of degree $n$ and $g:N\to L$ of degree $m$ we have
\[ d(g\circ f) \ = \ dg \circ f + (-1)^{m} g \circ df \]
as graded maps from $M$ to $L$.

The following consequences are crucial for the remaining step 
in the tilting theorem:\\
$\bullet$ for every DG $A$-module $M$, the endomorphism
complex $\End_A(M)=\Hom_A(M,M)$ is a differential graded
ring under composition and $M$ is a differential graded
$\End_A(M)$-$A$-bimodule;\\
$\bullet$ for every DG $A$-module $N$, the homomorphism
complex $\Hom_A(M,N)$ is a differential graded module
over $\End_A(M)$ under composition.
Moreover the functor $\Hom_A(M,-):\Mod A\to \Mod \End_A(M)$
is right adjoint to tensoring with the 
$\End_A(M)$-$A$-bimodule $M$;\\
$\bullet$ if $M$ is cofibrant, then the functor $\Hom_A(M,-)$
is exact and takes quasi-isomorphisms of DG $A$-modules
to quasi-isomorphisms. Moreover, its
left adjoint  $-\tensor_{\End_A(M)} M$ preserves cofibrant objects 
and chain homotopies.
So there exists a derived functor on the level of derived categories
\[ -\tensor^L_{\End_A(M)} M \ : \ \D(\End_A(M)) \ \to \  \D(A) \ , \] 
an exact functor which preserves infinite sums.

The following theorem is a special case of Lemma 6.1 in 
\cite{keller-derivingDG}.

\begin{theorem} \label{main theorem DG version} 
Let $A$ be a DG ring and $M$ a cofibrant $A$-module
which is a small generator for the derived category $\D(A)$.
Then the derived functor
\begin{equation}\label{eq-derived tensor with M}
 -\tensor^L_{\End_A(M)}M \ : \  \D(\End_A(M)) \ \to \ \D(A) \end{equation}
is an equivalence of triangulated categories.
\end{theorem}
\begin{proof}
The total left derived functor \eqref{eq-derived tensor with M}
is an exact functor between triangulated categories which
preserves infinite sums. Moreover, it takes the free 
$\End_A(M)$-module of rank one --- which is a small generator
for the derived category of $\End_A(M)$ --- to the small generator
$M$ for $\D(A)$. The induced map of graded endomorphism rings
\[ -\tensor^L_{\End_A(M)}M \, : \, \D(\End_A(M))(\End_A(M),\End_A(M))_*\ 
\to \ \D(A)(M,M)_* \]
is an isomorphism (both sides are isomorphic to the homology ring
of $\End_A(M)$). So Proposition \ref{triangulated comparison prop}
shows that this derived functor is an equivalence of 
triangulated categories.
\end{proof}

After all these preparations we can give the

\begin{proof}[Proof of the tilting theorem \ref{thm-tilting}]
Clearly, condition (3) implies condition (1).
Now we assume condition (1) and we choose an exact equivalence 
$F$ from the derived category $\D(R)$ to $\D(S)$.
The defining properties of a tilting complex are preserved under 
exact equivalences of triangulated categories.
Since $R[0]$, the free $R$-module of rank one, concentrated in dimension zero,
is a tilting complex for the ring $R$, its image  $T=F(R[0])$ 
is a tilting complex for $S$. Moreover, $F$ restricts to a ring isomorphism
\[ F \ : \ R \ \iso \ \D(R)(R[0],R[0]) \ \xrightarrow{\quad \iso\quad} \
 \D(S)(T,T) \ . \]
Hence condition (2) holds.

For the implication (2)$\Longrightarrow$(1) we are given a tilting complex  
$T$ in $\D(S)$ and an isomorphism of rings $\D(S)(T,T)\iso R$. 
The complex $T$ is naturally a differential graded $\End_S(T)$-$S$-bimodule,
and by Theorem \ref{main theorem DG version}, the derived functor
\[ -\tensor^L_{\End_S(T)} T \ : \  \D(\End_S(T)) \ \to \  \D(S)  \] 
is an equivalence of triangulated categories.
The isomorphism of graded rings
\[ H_* \left( \End_S(T)\right) \ \iso \ \D(S)(T,T)_* \] 
and the defining property of a tilting complex
show that the homology of $\End_S(T)$ is concentrated in dimension zero.
So there is a chain of two quasi-isomorphisms between 
$\End_S(T)$ and the ring 
$H_0\left(\End_S(T) \right) \iso \D(S)(T,T)\iso R$.
Restriction and extension of scalars along these quasi-isomorphisms
gives a chain of equivalences between the derived categories
of the differential graded ring $\End_S(T)$ and the derived category
of the ordinary ring $\D(S)(T,T)$.
Putting all of this together we end up with a chain of three equivalences
of triangulated categories:
\[ \D(R) \ \iso \ \D(\End_S(T)_+)\ \iso \ 
\D(\End_S(T)) \ \iso \  \D(S) \ .  \]

It remains to prove the implication (2)$\Longrightarrow$(3),
assuming that $R$ or $S$ is flat.
Let $T$ be a tilting complex in $\D(S)$ and $f:\D(S)(T,T)\to R$ 
an isomorphism of rings.
The homology of $\End_S(T)$ is isomorphic to
the graded self maps of $T$ in $\D(S)$, so it is concentrated in
dimension 0. So the inclusion of the DG sub-ring  $\End_S(T)_+$
into the endomorphism DG ring $\End_S(T)$ induces an isomorphism 
on homology, compare Remark \ref{rk-DG module remarks} (iv). 
Since $\End_S(T)_+$ is trivial in negative dimensions, 
there is a unique morphisms
of DG rings $\End_S(T)_+\to R$ which realizes the isomorphism $f$ on $H_0$.
We choose a {\em flat resolution} of  $\End_S(T)_+$, i.e., a DG ring $E$ 
and a quasi-isomorphism of DG rings $E\xrightarrow{\simeq} \End_S(T)_+$,
such that the functor $E\tensor -$ preserves quasi-isomorphisms between
chain complexes of abelian groups (see for example 
\cite[3.2 Lemma (a)]{keller-cyclic of exact}).
We end up with a chain of two quasi-isomorphisms of DG rings
\[ R \ \xleftarrow{\ \simeq\ } \ E \ \xrightarrow{\ \simeq\ } \
\End_S(T) \ .  \]
The complex $T$ is naturally a DG $\End_S(T)$-$S$-bimodule,
and we restrict the left action to  $E$ and view $T$ as a 
DG $E$-$S$-bimodule.
We choose  a cofibrant replacement $T^{\text{c}}\stackrel{\sim}{\to}T$ 
as a DG $E$-$S$-bimodule. Then we obtain the desired complex
of $R$-$S$-bimodules by
\[ M \ = \ R\tensor_E T^{\text{c}} \ . \] 
Tensoring with $M$ over $R$ has a total left derived functor
\begin{equation}\label{eq-derived tensor M}
 -\tensor_R^L M \ : \ \D(R) \ \to \ \D(S)  \end{equation}
(although this is not obvious with our definition since $M$ need not
be cofibrant as a complex of right $S$-modules, and then $-\tensor_RM$
does not takes values in cofibrant complexes).
In order to show that this derived functor is an exact equivalence
we use that the diagram of triangulated categories
\begin{equation}\label{eq-derived cat square}
\xymatrix@C=30mm@R=10mm{ \D(E) \ar[r]^-{-\tensor^L_E \End_S(T)}  
\ar[d]_{-\tensor^L_E R}&
\D(\End_S(T)) \ar[d]^{-\tensor^L_{\End_S(T)} T} \\
\D(R) \ar[r]_{ -\tensor_R^L M} & \D(S) } \end{equation}
commutes up to natural isomorphism.
Indeed, two ways around the square are given by derived tensor product
with the $E$-$S$-bimodules $M$ respectively $T$, so it suffices to
find a chain of quasi-isomorphisms of DG bimodules between $M$ and $T$.

Since $E$ is cofibrant as a complex of abelian groups and the composite map
$E\to \End_S(T)_+\to R$ is a quasi-isomorphism, 
$E\tensor S^{op}$ models the derived tensor product of $R$ and $S$. 
If one of $R$ or $S$ are flat, then $R\tensor S^{op}$ 
also models the derived tensor product, so that the map
\[ E\tensor S^{op} \ \to \ R\tensor S^{op} \] 
is a quasi-isomorphism of DG rings. Since $T^{\text{c}}$ is cofibrant as an 
$E\tensor S^{op}$-module, the induced map
\[ T^{\text{c}} \ = \  (E\tensor S^{op})_{E\tensor S^{op}} T^{\text{c}} 
\ \to \ 
 (R\tensor S^{op})_{E\tensor S^{op}} T^{\text{c}} \ \iso \
R\tensor_{E} T^{\text{c}} \ = \ M \]
is a quasi-isomorphism.
So we have a chain of two quasi-isomorphisms of $E$-$S$-bimodules
\[  T \ \xleftarrow{\ \simeq\ } \ T^{\text{c}} \ \xrightarrow{\ \simeq\ }  \ 
 M \ .  \] 

The left and upper functors in the commutative square
\eqref{eq-derived cat square} are derived from extensions of scalars
along quasi-isomorphisms of DG rings; thus they
are exact equivalence of triangulated categories.
The right vertical derived functor is an exact equivalence
by Theorem \ref{main theorem DG version}. 
So we conclude that the lower horizontal functor \eqref{eq-derived tensor M}
in the square \eqref{eq-derived cat square}
is also an exact equivalence of triangulated categories.
This establishes condition (3).
\end{proof}

\subsection{Examples} 
\label{sec-tilting examples}

Historically, tilting modules seem to have been the 
first examples of derived equivalences which are not Morita equivalences.
I will not try to give an account of the history of tilting module,
tilting complexes, and rather refer to \cite[Sec.~3.1]{koenig-zimmermann} or 
\cite{Morita obituary}.

\begin{defn}\label{def-tilting module}
Let $R$ be a finite dimensional algebra over a field. A {\em tilting module}
is a finitely generated $R$-module $T$ with the following properties.
\begin{itemize}
\item[(i)] $T$ has projective dimension 0 or 1,
\item[(ii)] $T$ has no self-extensions, i.e., $\Ext^1_R(T,T)=0$,
\item[(iii)] there is an exact sequence of right $R$-modules
\[ 0 \to R \to T_1 \to  T_2 \to 0 \] 
for some $m\geq 0$, such that $T_1$ and $T_2$ are direct summands of a finite
sum of copies of $T$.
\end{itemize}
\end{defn}

Note that if the tilting module $T$ is actually projective, then
condition (ii) is automatic and the exact sequence required in (iii)
splits. So then the free $R$-module of rank one is a summand of
a finite sum of copies of $T$, and hence $T$ is a finitely generated
projective generator for $\Mod R$. So $R$ is then Morita equivalent 
to the endomorphism ring of the tilting module $T$,
by the Morita theorem \ref{thm-classical Morita}.
For a self-injective algebra, for example a group algebra over a field,
the converse also holds;
indeed, every module of finite projective dimension over 
a self-injective algebra is already projective. So for these algebras,
tilting is the same as Morita equivalence.

If the projective dimension of the tilting module $T$ is 1, then we do not
get an equivalence between the modules over $R$ and $S=\End_R(T)$,
but we get a derived equivalence.
Indeed, since $R$ is noetherian, condition (i) implies that
$T$ has a 2-step resolution $P_1\to P_0$ by two finitely generated 
projective $R$-modules; this resolution is a small object in the derived
category $\D(R)$, and its graded self-maps in $\D(R)$ are concentrated in
dimension 0 by condition (ii).
Condition (iii) implies that the complex $R[0]$ is contained in the
triangulated subcategory generated by the resolution, and the resolution
is thus a generator for $\D(R)$, hence a tilting complex.

\begin{eg}\label{ex-tilting module A3 quivers}
For an example of a non-projective tilting module
we fix a field $k$ and we let $A$ be the algebra 
of upper triangular $3\times 3$ matrices over $k$,
\[ A \ =  \ \left\lbrace 
\begin{pmatrix}   x_{11} & x_{12} & x_{13} \\
                  0  & x_{22} & x_{23} \\
                  0  &   0  & x_{33} \end{pmatrix} \ | \
x_{ij} \in k \right\rbrace  \ . \]
Up to isomorphism, there are three indecomposable projective right $A$-modules,
namely the row vectors 
\[ P^1= \{ (y_{1}, y_2, y_3) \ | \ y_1,y_2,y_3\in k \} \]
and its $A$-submodules
\[ P^2= \{ (0, y_2, y_3) \ | \ y_2,y_3\in k \} 
\text{\quad and\quad}
P^3=\{ (0, 0, y_3) \ | \ y_3\in k \} \ . \]
These projectives are the covers of three corresponding simple modules,
namely 
\[ S^1= P^1/P^2 \ , \quad S^2= P^2/P^3 \ , \quad \text{and \quad} S^3=P^3\ .\] 
In particular, $S^3$ is projective and $S^1$ and $S^2$ have
projective dimension 1.

We define the tilting module $T$ as the direct sum
\[ T \ = \ P^1 \, \oplus \, P^2 \, \oplus \, S^2 \ . \] 
The projective resolution 
\[ 0 \ \xrightarrow{\qquad} \ P^1 \oplus P^2 \oplus P^3 \ 
\xrightarrow{\text{inclusion}} \ P^1 \oplus P^2 \oplus P^2 
\ \xrightarrow{\qquad} \ S^2 \ \xrightarrow{\qquad} \ 0 \]
can be used to calculate  $\Ext^1_A(T,T)=\Ext^1_A(S^2,T)=0$.
Since $P^1 \oplus P^2 \oplus P^3$ is a free $A$-module of rank one,
this short exact sequence verifies tilting condition (iii) in 
Definition \ref{def-tilting module} for the $A$-module $T$.

Altogether this shows that $T$ is a tilting module for $A$
of projective dimension one.
So $A$ `tilts' to the endomorphism algebra of $T$; 
this endomorphism algebra can be calculated directly,
and it comes out to be another subalgebra of the $3\times 3$ matrices 
over $k$, namely
\[ \End_A(T) \ \iso  \ \left\lbrace 
\begin{pmatrix}   x_{11} & x_{12} & x_{13} \\
                  0  & x_{22} & 0 \\
                  0  &   0  & x_{33} \end{pmatrix} \ | \
x_{ij} \in k \right\rbrace  \ . \]
(hint: the modules $P^1$ and $S^2$ do not map to each other nor to $P^2$, 
and the remaining relevant morphism spaces are 1-dimensional over $k$.)
The algebras $A$ and $\End_A(T)$ are {\em not} Morita equivalent.
Indeed, both have exactly three isomorphism classes of indecomposable
projective modules, but in one case these modules are `directed' 
(i.e., linearly ordered under the existence of non-trivial homomorphisms),
whereas in the other case two of these indecomposable
projectives do not map to each other non-trivially.

The preceding example, and many other ones, are often described 
using representations of quivers.
Indeed, the upper triangular matrices $A$ and the tilted algebra $\End_A(T)$
are isomorphic to the path algebras of the $A_3$-quivers
\[ \xymatrix{ \bullet \ar[r] & \bullet \ar[r] & \bullet } 
 \ \text{\qquad respectively\qquad}
\xymatrix{ \bullet & \bullet \ar[l]  \ar[r] & \bullet }  \ . \]
\end{eg}

\begin{eg}\label{ex-derived spread out matrix}
We obtain a tilting complex whose homology is concentrated 
in more then one dimension by
`spreading out' the free module of rank one.
Let $R=R_1\times R_2$ be the product of two rings.
Let $P_1=R_1\times 0$ and $P_2=0\times R_2$ be the two
"blocks", i.e., the projective $R$-bimodules corresponding to the central
idempotents $(1,0)$ and $(0,1)$ in $R$.
Then $R=P_1 \oplus P_2$ as an $R$-bimodule, and 
there are no non-trivial $R$-homomorphisms between $P_1$ and $P_2$.
Now take $T = P_1[0] \oplus P_2[n]$ for some number $n\ne 0$.
This is a complex of $R$-modules with trivial differential 
whose homology is concentrated in two dimensions. 
Moreover, the complex $T$ is a small generator for 
the derived category $\D(R)$. 
But the only non-trivial self-maps of $T$ are
of degree 0 since $P_1$ and $P_2 $ don't map to each other. 
Hence $T$ is a tilting complex which is not (quasi-isomorphic to) a 
tilting module. The endomorphisms of $T$ are again the ring $R$, 
so it is a non-trivial self-tilting complex of $R$.
Under the equivalence $\D(R)\iso D(R_1)\times D(R_2)$,
the self-equivalence induced by $T$ is the identity on
the first factor and the $n$-fold shift on the second factor.
\end{eg}

\bigskip

More examples of tilting complexes can be found in 
Sections 4 and 5 of~\cite{derived vs stable} or Chapter 5 
of~\cite{koenig-zimmermann}

\section{Morita theory in stable model categories}
\label{sec-Morita stable model categories}

Now we carry the Morita philosophy one step further:
we sketch Morita theory for ring spectra and for stable model categories.
As a summary one can say that essentially,
everything which we have said for rings and differential graded rings works, 
suitably interpreted, for ring spectra as well.

First a few words about what we mean by a ring spectrum.
The stable homotopy category of algebraic topology has a symmetric monoidal
smash product; the monoids are homotopy-associative ring spectra, and
they represent multiplicative cohomology theories. While the notion
of a homotopy-associative ring spectrum is useful for many things,
it does not have a good enough module theory for our present purpose.
One can certainly consider spectra with a 
homotopy-associative action of a homotopy-associative ring spectrum;
but the mapping cone of a homomorphism between such modules does not
inherit a {\em natural} action of the ring spectrum, 
and the category of such modules does not form a triangulated category.

So in order to carry out the Morita-theory program we need 
a highly structured model for the category of spectra which admits a symmetric
monoidal and homotopically well behaved smash product ---
before passing to the homotopy category ! 
The first examples of such categories were the $S$-modules~\cite{ekmm} 
and the symmetric spectra~\cite{hss};
by now several more such categories have been constructed
\cite{lydakis-simplicial, mmss};
The appropriate notion of a model category equivalence 
is a Quillen equivalence~\cite[Def.\ 1.3.12]{hovey-book} 
since these equivalences preserve the `homotopy theory',
not just the homotopy category; all known model categories
of spectra are Quillen equivalent in a monoidal fashion.

For definiteness, we work in one specific category of spectra 
with nice smash product, namely the {\em symmetric spectra} 
based on simplicial sets, 
as introduced by Hovey, Shipley and Smith~\cite{hss}.
The monoids are called {\em symmetric ring spectra}, and
I personally think that they are the simplest kind of ring spectra; 
as far as their homotopy category is concerned,
symmetric ring spectra are equivalent to the older notion of 
{\em $A_{\infty}$-ring spectrum}, 
and {\em commutative} symmetric ring spectra are equivalent 
to {\em $E_{\infty}$-ring spectra}. 
The good thing is that operads are not needed anymore.

However, using symmetric spectra is not essential and
the results described in this section could also be developed in 
more or less the same way in any other of the known 
model categories of spectra with compatible smash product.
Alternatively, we could have taken an axiomatic approach 
and use the term `spectra' for any stable, 
monoidal model category in which the unit object
`looks and feels' like the sphere spectrum. 
Indeed, we are essentially only using the
following properties of the category of symmetric spectra:
 
\begin{enumerate}
\item[(i)] there is a symmetric monoidal
{\em smash product}, which makes symmetric spectra into a 
{\em monoidal model category} (\cite[4.2.6]{hovey-book},
\cite{ss-monoidal});
\item[(ii)] the model structure is {\em stable}
(Definition \ref{def-stable model category});
\item[(iii)] the unit $\mS$ of the smash product is a small generator
(Definition \ref{def- triangulated small/generator})
of the homotopy category of spectra;
\item[(iv)] the (derived) space of self maps of the unit object $\mS$ is
weakly equivalent to $QS^0=\text{hocolim}_n\, \Omega^nS^n$ and
in the homotopy category, there are no maps of negative degree 
from $\mS$ to itself.
\end{enumerate}

A large part of the material in this section is taken from a
joint paper with Shipley~\cite{ss-classification}.
Two other papers devoted to Morita theory 
in the context of ring spectra are ~\cite{DGI} by 
Dwyer, Greenlees and Iyengar and \cite{BL} by Baker and Lazarev.

\subsection{Stable model categories}
\label{stable model categories}

Recall from ~\cite[I.2]{Q} or \cite[6.1]{hovey-book}
that the homotopy category of a pointed
model category supports a suspension and a loop functor. 
In short, for any object $X$ the map to the zero object can be factored
\[ X \ \xrightarrow{\quad } \ C \ \xrightarrow{\ \simeq\ } \ * \] 
as a cofibration followed by a weak equivalence. The suspension
of $X$ is then defined as the quotient of the cofibration,
$\Sigma X=C/X$. Dually, the loop object $\Omega X$ is the fiber
of a fibration from a weakly contractible object to $X$.
On the level of homotopy categories, the suspension and loop
constructions become functorial, and $\Sigma$ is left adjoint to $\Omega$.

\begin{defn} \label{def-stable model category}
A {\em stable model category} is a pointed model category for which
the functors $\Omega$ and $\Sigma$ on the homotopy category are
inverse equivalences.
\end{defn}

The homotopy category of a stable model category has  
a large amount of extra structure, some of which is relevant for us.
First of all, it is naturally a triangulated category,
see \cite[7.1.6]{hovey-book} for a detailed proof. 
The rough outline is as follows: by definition of `stable' the suspension 
functor is a self-equivalence of the
homotopy category and it defines the shift functor. 
Since every object is a two-fold suspension,
hence an abelian co-group object, the homotopy category of 
a stable model category is additive.
Furthermore, by \cite[7.1.11]{hovey-book} the cofiber sequences and 
fiber sequences of \cite[1.3]{Q} coincide
up to sign in the stable case, and they define the distinguished triangles. 
The model categories which we consider have all limits and colimits,
so the homotopy categories have infinite sums and products. 
Objects of a stable model category are called `generators' or `small' 
if they have this property as objects of the triangulated 
homotopy category, compare Definition \ref{def- triangulated small/generator}.

A Quillen adjoint functor pair between stable model categories 
gives rise to total derived functors which are exact functors
with respect to the triangulated structure; in other words 
both total derived functors commute with suspension and preserve
distinguished triangles. 

\begin{egs} \label{ex-stable model category examples}
\quad \\ (1)  {\bf Chain complexes.} In the previous section,
we have already seen an important class of examples from algebra, 
namely the category of chain complexes over a ring $R$. This category
actually has several different stable model structures:
the {\em projective} model structure 
(see Remark \ref{rk-cofibrant remark} (iv)) and the 
{\em injective} model structure (see Remark \ref{rk-cofibrant remark} (v))
have as weak equivalences the quasi-isomorphisms.
There is a clash of terminology here: the homotopy category 
in the sense of homotopical algebra is obtained by formally inverting the weak equivalences; so for the projective and  injective model structures, this 
gives the unbounded derived category $\D(R)$.
But the category of unbounded chain complexes admits another model
structure in which the weak equivalences are the {\em chain homotopy
equivalences}, see e.g.~\cite[Ex.~3.4]{Christensen-Hovey-relative}. 
Thus for this model structure, the homotopy category is what is 
commonly called the homotopy category, often denoted by $K(R)$.
The derived category $\D(R)$ is a quotient of the homotopy category $K(R)$; 
the derived category $\D(R)$ has a single small generator,
but for example the homotopy category of chain complexes of abelian groups  
$K(\mZ)$ does not have a set of generators whatsoever,
compare~\cite[E.3.2]{neeman-triangle book}. 
The three stable model structure on chain complexes of modules have
been generalized in various directions to chain complexes
in abelian categories or to other differential graded objects,
see ~\cite{Christensen-Hovey-relative}, \cite{beke-sheafy} 
and \cite{hovey-sheaves}.

(2)  {\bf The stable module category of a  Frobenius ring.} 
A different kind of algebraic example --- not involving chain complexes ---
is formed by the stable module categories of Frobenius rings. 
A Frobenius ring $A$ is defined by the property that  
the classes of projective and injective $A$-modules coincide. Important  
examples are  finite dimensional self-injective algebras over a field,
in particular finite dimensional Hopf-algebras, 
such as group algebras of finite groups.
The {\em stable module category} has as objects the $A$-modules 
({\em not} chain complexes of modules).
Morphisms in the stable category or represented by module homomorphisms,
but two homomorphisms are identified if their difference 
factors through a projective (= injective) $A$-module.

Fortunately the two different meanings of `stable' fit together nicely; 
the stable {\em module} category is the homotopy category
associated to a stable {\em model} category structure on the category of
$A$-modules, see~\cite[Sec.\ 2]{hovey-book}.  
The cofibrations are the monomorphisms, the fibrations are the epimorphisms,
and the weak equivalences are the maps which become isomorphisms in the
stable category.
Every finitely generated module is small when considered as an object 
of the stable module category.
As in the case of chain complexes of modules, there is usually no point in 
making the model structure explicit since the cofibration, fibrations and
weak equivalences coincide with certain well-known concepts.

Quillen equivalences between stable module categories 
arise under the name of {\em stable equivalences of Morita type}
(\cite[Sec.\ 5]{broue-block equivalences}, \cite[Ch.\ 11]{koenig-zimmermann}).
For simplicity, suppose that $A$ and $B$ are two finite-dimensional 
self-injective algebras over a field $k$; then $A$ and $B$ are in particular
Frobenius rings.
Consider an $A$-$B$-bimodule $M$ (by which we mean a $k$-symmetric
bimodule, also known as a right module over $A^{op}\tensor_k B$),
which is projective as left $A$-module and as a right $B$-module separately.
Then the adjoint functor pair
\begin{equation}\label{eq-Morita stable type pair} 
\xymatrix@=20mm{ \Mod A \quad \ar@<.4ex>^-{-\tensor_A M}[r] & 
\quad  \Mod B \quad  \ar@<.4ex>^-{\Hom_{B}(M,-)}[l] 
} \end{equation} 
is Quillen adjoint pair with respect to the `stable' model structures.

A stable equivalence of Morita type consists of an 
$A$-$B$-bimodule $M$ and a $B$-$A$-bimodule $N$ such that 
both $M$ and $N$ are projective as left and right modules
separately,
and such that there are direct sum decompositions
\[ N\tensor_A M \ \iso \ B\oplus X \text{\quad and \quad}
M\tensor_B N \ \iso \ A\oplus Y \] 
as bimodules, where $Y$ is a projective $A$-$A$-bimodule and
$B$ is a projective $B$-$B$-bimodule.
In this situation, the functors $-\tensor_A M$ and $-\tensor_BN$
induce inverse equivalences of the stable module categories.
Moreover, the Quillen adjoint pair \eqref{eq-Morita stable type pair} 
is a Quillen equivalence.

Rickard observed~\cite{derived vs stable} that a derived equivalence 
between self-injective, finite-dimen\-sional algebras
also gives rise to a stable equivalence of Morita type.

\medskip

In the above algebraic examples, there is no real need 
for the language of model categories;
moreover, `Morita theory' is covered by 
Keller's paper \cite{keller-derivingDG},
which uses differential graded categories.
A whole new world of stable model categories 
comes from homotopy theory, see the following list.
The associated homotopy categories yield 
triangulated categories which are not immediately visible 
to the eyes of an algebraist, since they do not arise from 
abelian categories.

\smallskip

(3)  {\bf Spectra.} 
The prototypical example of a stable model category (which is not
an `algebraic'), is `the' category of spectra. 
We review one model, the {\em symmetric spectra}
of Hovey, Shipley and Smith \cite{hss} in more detail in Section 
\ref{sec-symmetric spectra}.
Many other model categories of spectra have  
been constructed, see for example~\cite{BF, robinson-spectra, jardine-etale, 
ekmm, lydakis-simplicial, mmss}.
All known model categories of spectra Quillen equivalent
(see e.g.,  \cite[Thm.\ 4.2.5]{hss}, \cite{sch-comparison} or \cite{mmss}),
and their common homotopy category is referred to as {\em the}
stable homotopy category. The sphere spectrum is a small generator
for stable homotopy category.

(4)  {\bf Modules over ring spectra.} Modules over an  
$S$-algebra~\cite[VII.1]{ekmm}, over a symmetric ring spectrum  
\cite[5.4.2]{hss}, or over an orthogonal ring spectrum~\cite{mmss}
form stable model categories.
We recall symmetric ring spectra and their module spectra in 
Section \ref{sec-symmetric spectra}.
In each case a module is small if and only if it is weakly equivalent 
to a retract of a finite cell module. 
The free module of rank one is a small generator. 
More generally there are stable model categories of modules over 
`symmetric ring spectra with several objects', or {\em spectral categories}, 
see \cite[A.1]{ss-classification}.

(5)  {\bf Equivariant stable homotopy theory.} If $G$ is a compact  
Lie group, there is a category of $G$-equivariant coordinate free  
spectra~\cite{lms} which is a stable model category. 
Modern versions of this model category are the $G$-equivariant orthogonal
spectra of~\cite{mm} and $G$-equivariant $S$-modules of~\cite{ekmm}.
In this case  the equivariant suspension spectra of the coset spaces $G/H_+$ 
for all closed subgroups $H\subseteq G$ form a set of small generators.

(6)  {\bf Presheaves of spectra.}  For every Grothendieck site Jardine 
\cite{jardine-stable_presheaves}
constructs a stable model category of presheaves of 
Bousfield-Friedlander type spectra; the  weak equivalences are the maps 
which induce isomorphisms of the associated sheaves of stable homotopy groups.
For a general site these stable model categories do not seem to have a 
set of small generators. A similar model structure for
presheaves of symmetric spectra is developed 
in~\cite{jardine-symmetric_presheaves}.

(7)  {\bf The stabilization of a model category.}
Modulo technicalities, every pointed model category
gives rise to an associated stable model category by `inverting' 
the suspension functor, i.e., by passage to internal spectra.
This has been carried out, under different hypotheses, 
in \cite{sch-cotangent} and \cite{hovey-stabilization}.

(8)  {\bf Bousfield localization.} 
Following Bousfield~\cite{bousfield-localization-spaces}, localized model 
structures for modules over an $S$-algebra are constructed in 
\cite[VIII 1.1]{ekmm}. Hirschhorn \cite{hirschhorn-book}
shows that under quite general hypotheses the localization 
of a model category is again a model category. 
The localization of a stable model category is stable and localization
preserves generators. Smallness need not be preserved.

(9)  {\bf Motivic stable homotopy.} 
In \cite{morel-voevodsky, voevodsky-icm} Morel and Voevodsky introduced the 
${\mathbb A}^1$-local model category structure for schemes over a base. 
An associated stable homotopy category of ${\mathbb A}^1$-local 
$T$-spectra (where $T = {\mathbb A}^1/({\mathbb A}^1-0)$ is the `Tate-sphere') 
is an important tool in Voevodsky's proof of the Milnor conjecture 
\cite{voevodsky-milcon}. There are several stable model categories
underlying this motivic stable homotopy category, see for example
\cite{jardine-A^1stable}, \cite{hovey-stabilization}, 
\cite{hu-motivic S-modules} or~\cite{DOR}.
\end{egs}

\subsection{Symmetric ring and module spectra}
\label{sec-symmetric spectra}

In this section we give a quick introduction to symmetric spectra
and symmetric ring and module spectra. I recommend reading the original,
self-contained paper by Hovey, Shipley and Smith~\cite{hss}.
At several points, our exposition differs from theirs,
for example, we let the spheres act from the right.

\begin{defn}\cite{hss}
A symmetric spectrum consists of the following data
\begin{itemize}
\item a sequence of pointed simplicial sets
$X_n$ for $n\geq 0$
\item  for each $n\geq 0$ a base-point preserving action of
the symmetric group $\Sigma_n$ on $X_n$ 
\item pointed maps $\alpha_{p,q}:X_p\sm S^q\to X_{p+q}$ for $p,q\geq 0$
which are $\Sigma_p\times\Sigma_q$-equivariant;
here $S^1=\Delta^1/\partial\Delta^1$, $S^q=(S^1)^{\sm q}$ 
and $\Sigma_q$ permutes the factors.
\end{itemize}
This data is subject to the following conditions:
\begin{itemize}
\item under the identification $X_n\iso X_n\sm S^0$, the
map $\alpha_{n,0}:X_n\sm S^0\to X_n$ is the identity,
\item for $p,q,r\geq 0$, the following square commutes 
\begin{equation}\label{eq-symmetric axiom}\xymatrix@C=20mm{ 
X_p\sm S^q \sm S^r \ar^{\alpha_{p,q}\,\sm\,\Id}[r]
\ar_-{\iso}[d] & 
X_{p+q}\sm S^r \ar^{\alpha_{p+q,r}}[d]\\ 
X_p\sm S^{q+r} \ar_{\alpha_{p,q+r}}[r] & \ X_{p+q+r} \ . 
}\end{equation}
\end{itemize}
A {\em morphism} $f:X\to Y$ of symmetric spectra consists of 
$\Sigma_n$-equivariant pointed maps $f_n:X_n\to Y_n$ for $n\geq 0$,
which are compatible with the structure maps in the sense that
$f_{p+q}\circ\alpha_{p,q}=\alpha_{p,q}\circ (f_q\sm \Id_{S^q})$ 
for all $p,q\geq 0$.
\end{defn}

The definition we have just given is somewhat redundant, and
Hovey, Shipley and Smith use a more economical definition
in~\cite[Def.\ 1.2.2]{hss}.
Indeed, the commuting square \eqref{eq-symmetric axiom},
shows that all action maps  $\alpha_{p,q}$ are given by composites
of the maps $\alpha_{p,1}: X_p\sm S^1\to X_{p+1}$ for varying $p$.

A first example is the symmetric {\em sphere spectrum} $\mS$ given by
$\mS_n= S^n$, where the symmetric group permutes the factors and 
$\alpha_{p.q}:S^p\sm S^q\to S^{p+q}$ is the canonical isomorphism.
More generally, every pointed simplicial set $K$ gives rise to
a {\em suspension spectrum} $\Sigma^{\infty}K$ via
\[ (\Sigma^{\infty} K)_n \ = \ K\sm S^n \ ; \]
then we have $\mS\iso \Sigma^{\infty}S^0$.

A symmetric spectrum is {\em cofibrant} if it has the left lifting property
for levelwise acyclic fibrations. More precisely,  $A$ is cofibrant
if the following holds: for every morphism $f:X\to Y$ of symmetric spectra 
such that $f_n:X_n\to Y_n$ is a weak equivalence and Kan fibration for all $n$,
and for every morphism $\iota:A\to X$, there exists a morphism 
$\bar \iota:A\to Y$
such that $\iota=f\bar \iota$. Suspension spectra are examples of cofibrant
symmetric spectra.
An equivalent definition uses the {\em latching space} $L_nA$, 
a simplicial set which roughly is the 
`stuff coming from dimensions below $n$'; 
see~\cite[5.2.1]{hss} for the precise definition. 
A symmetric spectrum $A$ is cofibrant if and only if for all $n$, the map
$L_nA\to A_n$ is injective and symmetric group $\Sigma_n$ is freely
on the complement of the image, see~\cite[Prop. 5.2.2]{hss}. 
An {\em $\Omega$-spectrum} is defined by the properties that 
each simplicial set $X_n$ is a Kan complex 
and all the maps $X_n\to \Omega (X_{n+1})$ adjoint to
$\alpha_{n,1}$ are weak homotopy equivalences.
The {\em stable homotopy category} has as objects the cofibrant
symmetric $\Omega$-spectra and as morphisms the homotopy classes
of morphisms of symmetric spectra.

Although we just gave a perfectly good definition of the
stable homotopy category, in order to work with it 
one needs an ambient model category structure. One such model
structure is the {\em stable model structure} of \cite[Thm. 3.4.4]{hss}.
A morphism of symmetric spectra is a  {\em stable equivalence} if it induces
isomorphisms on all {\em cohomology theories} represented by
(injective) $\Omega$-spectra, see \cite[Def,\ 3.1.3]{hss} 
for the precise statement. There is a notion of 
{\em cofibration} such that a symmetric spectrum $X$ is cofibrant in 
the above sense if and only if the map from the trivial symmetric
spectrum to $X$ is a cofibration. The $\Omega$-spectra then 
coincide with the stably fibrant symmetric spectra.
There are other model structure for symmetric spectra with the same class of
weak (=stable) equivalences, hence with the same homotopy category,
for example the {\em $S$-model structure} 
which is hinted at in~\cite[5.3.6]{hss}.

{\bf Stable equivalences versus  $\pi_*$-isomorphisms.} 
One of the tricky points with symmetric spectra is the relationship
between stable equivalences and $\pi_*$-isomorphisms.
The stable equivalences are defined as the morphisms which induce
isomorphisms on all {\em cohomology theories}; there is the
strictly smaller class of morphisms which induce isomorphisms on
stable homotopy groups.
The {\em $k$-th stable homotopy group} of a symmetric spectrum $X$ is
defined as the colimit
\[ \pi_k X \ = \ \colim_n \, \pi_{n+k} |X_n| \ , \] 
where $|X_n|$ denotes the geometric realization of the simplicial set $X_n$. 
The colimit is taken over the maps
\begin{equation}\label{colimit system pi_n}  
\pi_{n+k} \, |X_n| \ \xrightarrow{\ -\sm S^1\ } \
\pi_{n+k+1} \, \left( |X_{n}|\sm S^1 \right) \ \xrightarrow{\ (\alpha_{n,1})_*\ } \
\pi_{n+k+1} \, |X_{n+1}| \ . \end{equation}

While every $\pi_*$-isomorphism of symmetric spectra is a stable
equivalence~\cite[Thm.\ 3.1.11]{hss}, the converse is not true.
The standard example is the following:
consider the symmetric spectrum $F_1S^1$ freely generated by the circle $S^1$
in dimension 1. Explicitly,  $F_1S^1$ is given by
\[ (F_1S^1)_n \ = \ \Sigma_n^+\sm_{\Sigma_{n-1}}S^{n-1}\sm S^1 \ . \]
So $(F_1S^1)_n$ is a wedge of $n$ copies of $S^n$
and in the stable range, i.e., up to roughly dimensions $2n$, 
the homotopy groups of $(F_1S^1)_n$ are a direct sum of
$n$ copies of the homotopy groups of $S^n$. Moreover, in the stable range,
the map in the colimit system \eqref{colimit system pi_n} 
is a direct summand inclusion into $(n+1)$ copies of the homotopy groups 
of $S^n$. Thus in the colimit, the stable homotopy groups of the
symmetric spectrum  $F_1S^1$ are a countably infinite direct sum of
copies of the stable homotopy groups of spheres.
Since $F_1S^1$ is freely generated by the circle $S^1$ in dimension 1,
it ought to be a desuspension of the suspension spectrum of the circle. 
However, the necessary symmetric group actions `blow up'
such free objects with the effect that the stable
homotopy groups are larger than they should be.
This example indicates that inverting only the  $\pi_*$-isomorphisms
would leave too many stable homotopy types, and the resulting category
could not be equivalent to the usual stable homotopy category.

{\bf Smash product.}
One of the main features which distinguishes {\em symmetric} spectra 
from the more classical spectra is the internal smash product.
The smash product of symmetric spectra can be described
via its universal property, analogous to universal property 
of the tensor product over a commutative ring.
Indeed, if $R$ is a commutative ring and $M$ and $N$ are right $R$-modules,
then a bilinear map to another left $R$-module $W$ is a map $b:M\times N\to W$
such that for each $m\in M$ the map $b(m,-):N\to W$
and each $n\in N$ the map $b(-,n):M\to W$ are $R$-linear.
The tensor product $M\tensor_R N$ is the universal example of
a right $R$-module together with a bilinear map from $M\times N$.

Let us define a {\em bilinear morphism} $b:(X,Y)\to Z$ 
from two symmetric spectra $X$ and $Y$
to a symmetric spectrum $Z$ to consist of a collection of
$\Sigma_p\times\Sigma_q$-equivariant maps of pointed simplicial sets
\[ b_{p,q} \ : \ X_p \sm \ Y_q \ \to \ Z_{p+q} \]
for $p,q\geq 0$, such that
for all $p,q,r\geq 0$, the following diagram commutes
\begin{equation}\label{eq-bilinear diagram}
\xymatrix@C=15mm@R=10mm{ & 
X_p \sm Y_q \sm S^r \ar[dl]_{\Id\sm \alpha_{q,r}} 
\ar[d]^{b_{p,q}\sm \Id} 
\ar[r]^-{\Id\sm \text{twist}}& 
X_p \sm \ S^r \sm Y_q \ar[d]^{\alpha_{p,r}\sm\Id}\\
X_{p} \sm Y_{q+r}  \ar[dr]_{b_{p,q+r}} & 
Z_{p+q}\sm S^r \ar[d]^{\alpha_{p+q,r}} & X_{p+r}\sm Y_{q} \ar[d]^{b_{p+r,q}} \\
&  Z_{p+q+r}  & Z_{p+r+q}\ar[l]^-{1\times\chi_{r,q}\times 1}}\end{equation}
The automorphism $1\times\chi_{r,q}$
of $Z_{p+q+r}$ may look surprising at first sight. Here
$1\times \chi_{r,q}\in \Sigma_{p+r+q}$ denotes the block permutation
which fixes the first $p$ elements, and which moves the next $q$ 
elements past the last $r$ elements. 
This can be viewed as a topological version
of the Koszul sign rule which says that when two symbols of degree
$q$ and $r$ are permuted past each other, the sign $(-1)^{qr}$
should appear as well. The block permutation $\chi_{r,q}$ has sign
$(-1)^{qr}$ and it compensates the upper vertical interchange of
$Y_q$ and $S^r$.
A good way to keep track of such permutations is to carefully distinguish
between indices such as $r+q$ and $q+r$. Of course these two numbers are
equal, but the fact that one arises naturally instead of the other
reminds us that a block permutation should be inserted.

The smash product $X\sm Y$ is the universal example 
of a symmetric spectrum with a bimorphism from $X$ and $Y$.
In other words, it comes with a bimorphism $\iota:(X,Y)\to X\sm Y$
such that for every symmetric spectrum $Z$ the map
\begin{equation}\label{eq-universal property smash}
 \spec^{\Sigma}(X\sm Y,Z) \ \to \ \text{Bi-}\spec^{\Sigma}((X,Y),Z) 
\end{equation}
is bijective. If we suppose that such a universal object exist,
this property characterizes the smash product and the maps
$\iota_{p,q}:X_p\sm Y_q\to (X\sm Y)_{p+q}$ up to canonical isomorphism.
An actual construction as a certain coequalizer is given 
in~\cite[Def.\ 2.2.3]{hss}; in this article, we will only use 
the universal property of the smash product.

We use the universal property to derive that the smash product
is functorial and symmetric monoidal.
For example, let $f:X\to Y $ and $f':X'\to Y'$ be morphisms 
of symmetric spectra. Then the collection of maps of pointed simplicial sets
\[ \left\lbrace X_p \sm \ X'_q \ \xrightarrow{f_q\sm f'_q} \  
Y_p \sm \ Y'_q \ \xrightarrow{\iota_{p,q}} \ (Y\sm Y')_{p+q} 
\right\rbrace_{p,q\geq 0} \]
form a bilinear morphism $(X,X')\to Y\sm Y'$, so it corresponds
to a unique morphism of symmetric spectra $f\sm f':X\sm X'\to Y\sm Y'$.
The universal property implies functoriality in both arguments.
For the proof of the associativity of the smash product we notice that
 the family
\[ \left\lbrace X_p \sm \ Y_q \sm Z_r \ \xrightarrow{\iota_{p,q}\sm \Id} \  
(X\sm Y)_{p+q} \sm Z_r  \ \xrightarrow{\iota_{p+q,r}} \
 ((X\sm Y)\sm Z)_{p+q+r} \right\rbrace_{p,q,r\geq 0} \]
and the family
\[ \left\lbrace X_p \sm \ Y_q \sm Z_r \ \xrightarrow{\Id\sm\iota_{q,r}} \  
X_p\sm (Y\sm Z)_{q+r}  \ \xrightarrow{\iota_{p,q+r}} \
 (X\sm (Y\sm Z))_{p+q+r} \right\rbrace_{p,q,r\geq 0} \]
both have the universal property of a {\em tri-linear} morphism 
out of $X$, $Y$ and $Z$.
The uniqueness of universal objects gives a preferred isomorphism
of symmetric spectra
\[ (X\sm Y)\sm Z \ \iso X\sm (Y\sm Z)  \ . \] 
The symmetry isomorphism $X\sm Y\iso Y\sm X$ corresponds to the
bilinear morphism
\begin{equation}\label{eq-symmetry iso}
 \left\lbrace X_p \sm \ Y_q \ \xrightarrow{\text{twist}} \  
Y_q \sm \ X_p\ \xrightarrow{\iota_{q,p}} \ (Y\sm X)_{q+p} 
\xrightarrow{\chi_{q,p}} \ (Y\sm X)_{p+q} 
\right\rbrace_{p,q\geq 0} \ . \end{equation}
The block permutation $\chi_{q,p}$ is crucial here: without
it we would not get a bilinear morphism is the sense of diagram
\eqref{eq-bilinear diagram}.
In much the same spirit, the universal properties can be used to 
provide unit isomorphisms $\mS\sm X\iso X\iso X\sm \mS$, to verify the
coherence conditions of a symmetric monoidal structure, and to establish
an isomorphism of suspension spectra
\[ (\Sigma^{\infty}K)\, \sm\, (\Sigma^{\infty}L) \ \iso \ 
\Sigma^{\infty}(K\sm L) \ . \]
The symmetric monoidal structure given by the smash product of
symmetric spectra is {\em closed} in the sense that internal function
objects exist as well. For each pair of symmetric spectra $X$ and $Y$ there
is a symmetric function spectrum $\Hom(X,Y)$~\cite[2.2.9]{hss}, 
and there are natural composition morphisms
\[ \circ \ : \ \Hom(Y,Z) \, \sm \, \Hom(X,Y) \ \to \ \Hom(X,Z) \] 
which are associative and unital with respect to a unit map $\mS\to \Hom(X,X)$.
Moreover, the usual adjunction isomorphism
\[ \spec^{\Sigma}(X\sm Y, Z) \ \iso \ \spec^{\Sigma}(X,\Hom(Y,Z)) \] 
relates the smash product and function spectra.

{\bf Ring and module spectra.} The smash product of symmetric spectra leads to
the concomitant concepts  {\em symmetric ring spectra}, {\em module spectra}
and {\em algebra spectra}.

\begin{defn} \label{symmetric ring spectrum} 
A {\em symmetric ring spectrum} is a symmetric spectrum $R$ together 
with morphisms of symmetric spectra
\[  \eta \ \colon \ \mS\ \to \ R \quad \mbox{and} \quad 
\mu \ \colon \ R\,\sm\, R \ \to\ R \ , \]
called the unit and multiplication map, which satisfy certain 
associativity and unit conditions (see~\cite[VII.3]{MacLane-working}). 
A ring spectrum
$R$ is {\em commutative}  if the multiplication map is unchanged 
when composed with the twist, or the symmetry isomorphism
\eqref{eq-symmetry iso}, of $R\,\sm\, R$. 
A morphism of ring spectra is a morphism of spectra 
commuting with the multiplication and unit maps. 
If $R$ is a symmetric ring spectrum, a {\em right $R$-module} 
is a spectrum $N$ together with an action map $N\,\sm\, R\to N$ 
satisfying associativity and unit conditions 
(see again~\cite[VII.4]{MacLane-working}).
A morphism of right $R$-modules is a morphism of spectra commuting 
with the action of $R$. We denote the category of right $R$-modules 
by $\Mod R$.

\end{defn}

With the universal property of smash product we can make the structure
of a symmetric ring spectrum more explicit. The multiplication map
$\mu:R\sm R\to R$ corresponds to a family of pointed, 
$\Sigma_p\times\Sigma_q$-equivariant maps
\[ \mu_{p,q} \ : \ R_p \sm \ R_q \ \to \ R_{p+q} \]
for $p,q\geq 0$, which are bilinear in the sense of diagram 
\eqref{eq-bilinear diagram}. The maps are supposed to be associative
and unital with respect to the maps $\eta_p:S^p\to R_p$ which constitute the
unit map $\eta:\mS\to R$.

The commutativity isomorphism of the smash product involves 
the block permutation $\chi_{q,p}$, see \eqref{eq-symmetry iso}.
So the multiplication of a symmetric ring spectrum is commutative
if and only if the following diagrams commute for all $p,q\geq 0$
\[\xymatrix@C=15mm{  R_p \sm \ R_q \ar[r]^-{\mu_{p,q}} \ar[d]_{\text{twist}} &
R_{p+q} \ar[d]^{\chi_{p,q}} \\
R_q \sm \ R_p \ar[r]_-{\mu_{q,p}} & R_{q+p} } \]
The block permutation $\chi_{p,q}$ has sign $(-1)^{pq}$, so this
diagram is reminiscent of the Koszul sign rule in a graded ring
which is commutative in the graded sense.

\bigskip

The unit $\mS$ of the smash product is a ring spectrum in a unique way, 
and $\mS$-modules are the same as symmetric spectra. 
The smash product of two ring spectra is naturally a ring spectrum. 
For a ring spectrum $R$ the opposite ring spectrum $R^{\text{op}}$ is defined 
by composing the multiplication with the twist map $R\,\sm\, R\to R\,\sm\, R$
(so in terms of the bilinear maps $\mu_{p,q}:R_p\sm R_q\to R_{p+q}$,
a block permutation appears).
The definitions of left modules and bimodules is hopefully clear;
left $R$-modules and $R$-$T$-bimodule can also be defined as right modules 
over the opposite ring spectrum $R^{op}$, respectively right modules
over the ring spectrum $R^{\text{op}}\,\sm\, T$.

A formal consequence of having a closed symmetric monoidal smash product 
is that the category of $R$-modules inherits a smash product 
and function objects. The smash product $M\sm_R N$ of  a right $R$-module $M$
and a left $R$-module $N$ can be defined as the coequalizer, 
in the category of symmetric spectra, of the two maps 
\[ \xymatrix{M\,\sm\, R\,\sm\, N \ar@<.3ex>[r] \ar@<-.3ex>[r] & M\,\sm\, N }\]
given by the action of $R$ on $M$ and $N$ respectively. 
Alternatively, one can characterize  $M\sm_R N$ as the universal example
of a symmetric spectrum which receives a bilinear map from
$M$ and $N$ which is {\em $R$-balanced}, i.e., all the diagrams
\begin{equation}\label{eq- R-module bilinear}
\xymatrix@C=15mm{ 
M_p \sm R_q \sm N_r \ar[d]_{\alpha_{p,q}\sm \Id} \ar[r]^{\Id\sm\alpha_{q,r}} &
M_p\sm N_{q+r} \ar[d]^{\iota_{p,q+r}} \\
M_{p+q}\sm N_r \ar[r]_{\iota_{p+q,r}} & \ (M\sm N)_{p+q+r} }\end{equation}
commute. If $M$ happens to be a $T$-$R$-bimodule and $N$ an $R$-$S$-bimodule, 
then $M\sm_R N$ is naturally a $T$-$S$-bimodule. In particular, 
if $R$ is a commutative ring spectrum, the notions of left 
and right module coincide and agree with the notion of a symmetric bimodule. 
In this case $\sm_R$ is an internal symmetric monoidal smash product 
for $R$-modules. There are also internal function spectra and
function modules, enjoying the `usual' adjointness properties
with respect to the various smash products.

The modules over a symmetric  ring spectrum $R$
inherit a model category structure from symmetric spectra, 
see~\cite[Cor. 5.4.2]{hss} and~\cite[Thm.\ 4.1 (1)]{ss-monoidal}.
More precisely, a morphism of $R$-modules is 
called a weak equivalence (resp.\ fibration) if the underlying morphism
of symmetric spectra is a stable equivalence (resp.\ stable fibration).
The cofibrations are then determined by the left lifting property with
respect to all acyclic fibrations in $\Mod R$.
This model structure is stable, so the homotopy category of modules 
over a ring spectrum is a triangulated category. 
The free module of rank one is a small generator.

For a map $R\to S$ of ring spectra, there is a Quillen adjoint 
functor pair analogous to restriction and extension of scalars: 
any $S$-module becomes an $R$-module if we let $R$ act through the map. 
This functor has a left adjoint taking an $R$-module $M$ to the $S$-module 
$M\, \sm_R \, S$. If $R\to S$ is a stable equivalence, 
then the functors of restriction and extension of scalars are 
a Quillen equivalence between the categories of $R$-modules and $S$-modules,
 see~\cite[Thm.\ 5.4.5]{hss} 
and~\cite[Thm.\ 4.3]{ss-monoidal}.

\begin{eg}[Monoid ring spectra]
If $M$ is a simplicial monoid, and $R$ is a symmetric ring spectrum, 
we define a symmetric spectrum $R[M]$ by
\[ R[M]_n \ = \ R_n \sm M^+ \ , \]
where $M^+$ denotes the underlying simplicial set of $M$, 
with disjoint basepoint added.
The unit map is the composite 
of the unit map of $R$ and the wedge summand inclusion
indexed by the unit of $M$; the multiplication map 
$R[M]\sm R[M]\to R[M]$ is induced from the bilinear morphism
\[ (R_p\sm M^+) \sm  (R_q\sm M^+) \iso  (R_p\sm R_q)\sm (M\times M)^+ 
\ \xrightarrow{\ \mu_{p,q}\sm \text{mult.}\ } \ R_{p+q} \sm M^+ \ . \]
The construction of the monoid ring over $\mS$ is left adjoint 
to the functor which takes a symmetric ring spectrum $R$ to the 
simplicial monoid $R_0$.
\end{eg}

\begin{eg}[Matrix ring spectra] \label{ex-matrix ring spectrum}
Let $R$ be a symmetric ring spectrum and consider the wedge (coproduct)
\[ R\times n \ = \ 
\underbrace{R\vee \dots \vee R}_n \]
of $n$ copies of the free $R$-module of rank 1. 
In the usual stable model structure, the free module of rank 1 is 
cofibrant, hence so is $R\sm n^+$. We choose a fibrant replacement 
$R \times n\xrightarrow{\sim} (R \times n)^{\text{f}}$.
The ring spectrum of {\em $n\times n$ matrices} over $R$ is defined as the
endomorphism ring spectrum of $(R \times n)^{\text{f}}$, 
\[ M_n(R) \ = \ \End_{R}((R \times n)^{\text{f}}) \ . \] 
The stable equivalence type of the matrix ring spectrum $M_n(R)$
is independent of the choice of fibrant replacement,
see Corollary A.2.4 of \cite{ss-classification}.
Moreover, the underlying spectrum of $M_n(R)$ is isomorphic, in the
stable homotopy category, to a sum of $n^2$ copies of $R$.
\end{eg}

\begin{eg}[Eilenberg-Mac Lane spectra]
\label{ex-EM ring spectra}

For an abelian group $A$, the {\em Eilenberg-Mac Lane spectrum} $HA$
is defined by
\[ (HA)_n \ = \ A\tensor \mZ[S^n] \ , \] 
i.e., the underlying simplicial set of the
dimensionwise tensor product of $A$ with the reduced free simplicial abelian
generated by the simplicial $n$-sphere.
The symmetric groups acts by permuting the smash factors of $S^n$.
The homotopy groups of the symmetric spectrum $HA$ are concentrated
in dimension zero, where we have a natural isomorphism
$\pi_0HA\iso A$.

For two abelian groups $A$ and $B$, a natural morphism 
of symmetric spectra 
\[ HA\sm HB \ \to \  H(A\tensor B) \] 
is obtained, by the universal property \eqref{eq-universal property smash}, 
from the bilinear morphism
\begin{align*} (HA)_n \sm (HB)_m \ = \ 
\left( A\tensor \mZ[S^n]\right) \, &\sm \,\left( B\tensor\mZ[S^m]\right) \\
 &\to \ (A\tensor B)\tensor \mZ[S^{n+m}] \ = \ (H(A\tensor B))_{n+m} 
\end{align*}
given by
\[ \left( \sum_i a_i\cdot x_i \right) \, \sm \, 
\left( \sum_j b_j\cdot x'_j \right) \ \longmapsto \ 
\sum_{i,j} (a_i\cdot b_j)\cdot \, x_i\sm x'_j \  . \] 
A unit map $\mS\to H\mZ$ is given by the inclusion of generators.
With respect to these maps, $H$ becomes a 
lax symmetric monoidal functor from the
category of abelian groups to the category of symmetric spectra.
As a formal consequence, $H$ turns a ring $R$ into a symmetric ring spectrum
with multiplication map
\[ HR\sm HR \ \to \ H(R\tensor R) \ \to \ HR \ . \]
Similarly, an $R$-module structure on $A$ gives rise to
an $HR$-module structure on $HA$.

The definition of the symmetric spectra $HA$ makes just as much sense
when $A$ is a {\em simplicial} abelian group; 
thus the Eilenberg-Mac Lane functor makes simplicial rings into symmetric
ring spectra, respecting possible commutativity of the multiplications.
With a little bit of extra care, the Eilenberg-Mac Lane construction 
can also be extended to a 
differential graded context, compare~\cite[App. B]{ss-classification}
and \cite{richter-DoldKan}.

For a fixed ring $B$, the modules
over the Eilenberg-Mac Lane ring spectrum $HB$
of a ring $B$ have the same homotopy theory as complexes of $B$-modules.
The first results of this kind were obtained 
by Robinson for $A_{\infty}$-ring spectra~\cite{robinson-derived}, and
later for $S$-algebras in~\cite[IV Thm.\ 2.4]{ekmm}; in both cases, 
equivalences of triangulated homotopy categories are constructed. 
But more is true: for any ring $B$, 
Theorem 5.1.6 of \cite{ss-classification}
provides a chain of two Quillen equivalences between 
the categories of unbounded chain complexes of $B$-modules and
the $HB$-module spectra.
\end{eg}

\begin{eg}[Cobordism spectra] We define a commutative symmetric ring spectrum
$MO$ whose stable homotopy groups are isomorphic to the ring of cobordism
classes of closed manifolds.
We set
\[  (MO)_n \ = \ EO(n)^+\sm_{O(n)} S^n \ . \]
Here $O(n)$ is the $n$-th orthogonal group consisting of
Euclidean automorphisms of $\mR^n$.
The space $EO(n)$ is the geometric realization of the simplicial object
of topological groups which in dimension $k$ is the $k$-fold product
of copies of $O(n)$, and where are face maps are projections.
Thus  $EO(n)$ is a topological group with a homomorphism $O(n)\to EO(n)$
coming from the inclusion of 0-simplices. The underlying space of $EO(n)$
is contractible and has two commuting actions of $O(n)$ from the left and
the right. The right $O(n)$-action is used to form the orbit space
$(MO)_n$, where we think of $S^n$ as the one-point compactification
of $\mR^n$with its natural left  $O(n)$-action. 
Thus the space $(MO)_n$ still has a left  $O(n)$-action, which we
restrict to an action of the symmetric group $\Sigma_n$,
sitting inside $O(n)$ as the coordinate permutations.
Topologically, $(MO)_n$ is nothing but the Thom space of the tautological
bundle over the space $BO(n)$.

The unit of the ring spectrum $MO$ is given by the maps
\[ S^n \iso  O(n)^+\sm_{O(n)} S^n \ \to \ 
EO(n)^+\sm_{O(n)} S^n = (MO)_n \]  using the `vertex map' $O(n)\to EO(n)$.
There are multiplication  maps 
\[ (MO)_p \sm (MO)_q \ \to \ (MO)_{p+q} \] 
which are induced from the identification $S^p\sm S^q\iso S^{p+q}$
which is equivariant with respect to the group $O(p)\times O(q)$,
viewed as a subgroup of $O(p+q)$.
The fact that these multiplication maps are associative and commutative
uses that
\begin{itemize}
\item for topological groups $G$ and $H$, the simplicial model of $EG$ 
comes with a natural,
associative and commutative isomorphism $E(G\times H)\iso EG\times EH$;
\item the group monomorphisms  $O(p)\times O(q)\to O(p+q)$
are strictly associative, and the following diagram commutes
 \[\xymatrix@C=15mm{  O(p)\times O(q) \ar[r]^-{}
\ar[d]_{\text{twist}} &
O(p+q) \ar[d]^{\text{conj. by }\chi_{p,q}} \\
O(q) \times O(p) \ar[r]_-{} & O(q+p) } \]
where the right vertical map is conjugation by the permutation matrix
of the block permutation $\chi_{p,q}$.
\end{itemize}
In very much the same way we obtain commutative symmetric ring spectra
model for the oriented cobordism spectrum $MSO$ and the
spin cobordism spectrum $MSpin$. The complex cobordism ring spectrum 
$MU$ does not fit in here directly; one has to vary the notion of
a symmetric spectrum slightly, and consider only 
symmetric spectra which are defined `in even dimensions'.
\end{eg}

\subsection{Characterizing module categories over ring spectra}

Several of the examples of stable model categories in
Section \ref{stable model categories} already come as categories of modules 
over suitable rings or ring spectra. This is no coincidence.
In fact, every stable model category with a single small
generator has the same homotopy theory as the modules over a ring spectrum.
This is an analog of Theorem \ref{thm-Gabriel/Michell}, which
characterizes module categories over a ring as the 
cocomplete abelian category with a small projective generator.

To an object $P$ in a sufficiently nice stable model category $\C$ we
can associate a symmetric {\em endomorphism ring spectrum} $\End_{\C}(P)$;
among other things, this ring spectrum comes with an isomorphism 
of graded rings
\[ \pi_* \End_{\C}(P) \ \iso \ {\Ho(\C)}(P,P)_* \] 
between the homotopy groups of  $\End_{\C}(P)$ and the morphism of $P$ in the
triangulated homotopy category of $\C$.
The precise result is as follows:

\begin{theorem} \label{thm-main-one}
Let $\C$ be a  stable model category which is simplicial,
cofibrantly generated and proper.
If $\C$ has a small generator $P$, then there exists 
a chain of Quillen equivalences between $\C$ and the model category
of $\End_{\C}(P)$-modules,
\[ \C \simeq_Q \Mod\End_{\C}(P) \ . \]
\end{theorem}

Unfortunately, the theorem is currently only known 
under the above technical hypothesis: 
the stable model category in question should be
simplicial (see~\cite[II.2]{Q}, \cite[4.2.18]{hovey-book}),  
{\em cofibrantly generated} (see~\cite[Sec.\ 2.1]{hovey-book}
or~\cite[Sec.\ 2]{ss-monoidal}) and proper 
(see \cite[Def.\ 1.2]{BF} or \cite[Def.\ 5.5.2]{hss}).
The conditions enter in the construction of `good' endomorphism ring spectra.
I suspect however, that these hypothesis are not essential and
can be eliminated with a clever use of framing techniques.
For example, in \cite[Sec.\ 6]{ss-uniqueness},
we use framings to construct function spectra in a 
arbitrary stable model category; that construction does however
not yield {\em symmetric} spectra, and there is no good composition pairing.
The condition of being a {\em simplicial} model category can be removed
by appealing to \cite{rss} or \cite{dugger-simplicial} where
suitable model categories are replaced by Quillen equivalent simplicial
model categories.

This theorem is a special case of the more general result which applies to 
stable model categories with a set of small generators
(as opposed to a single small generator), 
see \cite[Thm.\ 3.3.3]{ss-classification}.

{\bf Spectral model categories.}
In the algebraic situations which we considered 
in Sections \ref{sec-classical Morita} and \ref{sec-derived Morita}, 
the key point is to have a good notion of endomorphism ring 
or endomorphism DG ring
together with a `tautological' functor
\begin{equation}\label{Hom(P,-)} 
\Hom_{\A}(P,-) \ : \ \A \ \to \ \Mod \End_{\A}(P) \ . \end{equation}
Then it is a matter of checking that when $P$ is a small
generator, the functor $\Hom(P,-)$ is either an equivalence of categories
(in the context of abelian categories) or induces an equivalence of
derived categories (in the context of DG categories).
For abelian categories the situation is straightforward, and the ordinary
endomorphism ring does the job. In the differential graded context already
a little complication comes in because the categorical hom functor
$\Hom_{A}(P,-)$ need not preserve quasi-isomorphisms in general.

For stable model categories, the key construction is again to have
an {\em endomorphism ring spectrum} $\End_{\C}(P)$ 
together with a homotopically
well-behaved homomorphism functor \eqref{Hom(P,-)} to modules
over the endomorphism ring spectrum.
This is easy for the following class of {\em spectral model categories} 
where composable function spectra are part of the data.
A spectral model category is analogous 
to a simplicial model category~\cite[II.2]{Q}, 
but with the category of simplicial sets replaced 
by symmetric spectra.  
Roughly speaking, a spectral model category is a pointed model category 
which is compatibly enriched over the stable model category of spectra. 
In particular there are `tensors' $K\sm X$ and `cotensors' $X^K$ 
of an object $X$ of $\C$ and a symmetric spectrum $K$, and 
function symmetric spectra $\Hom_{\C}(A,Y)$ between two objects of $\C$.
The compatibility is expressed by the following axiom
which takes the place of \cite[II.2 SM7]{Q}; 
there are two equivalent `adjoint' forms of this axiom,
compare \cite[Lemma 4.2.2]{hovey-book} or \cite[3.5]{ss-classification}.

(Pushout product axiom) 
For every cofibration $A\to B$ in $\C$ and every cofibration $K\to L$  
of symmetric spectra, the {\em pushout product map}
\[ L \sm A \cup_{K \sm A}  K \sm B \ \to \ L \sm B \]
is a cofibration; the pushout product map is a weak equivalence if  
in addition $A\to B$ is a weak equivalence in $\C$ or $K\to L$ is a  
stable equivalence of symmetric spectra.

A {\em spectral Quillen pair} is a Quillen adjoint functor pair 
$L:\C\to \D$ and $R:\D\to\C$ between spectral model categories
together with a natural isomorphism of symmetric homomorphism spectra
\[ \Hom_{\C}(A,RX) \ \iso \  \Hom_{\D}(LA,X) \]
which on the vertices of the 0-th level reduces to the adjunction isomorphism.
A spectral Quillen pair is a {\em spectral Quillen equivalence} if
the underlying Quillen functor pair is an ordinary Quillen equivalence.

A spectral model category is the  same as a `$\spec^\Sigma$-model category' 
in the sense of \cite[Def.\ 4.2.18]{hovey-book};
Hovey's condition 2 of \cite[4.2.18]{hovey-book} is
automatic since the unit $\mS$ for the smash product of symmetric spectra 
is cofibrant.
Similarly, a spectral Quillen pair is a `$\spec^\Sigma$-Quillen functor' 
in Hovey's terminology.
Examples of spectral model categories are module categories over
a ring spectrum, and the category of symmetric spectra over 
a suitable simplicial model category~\cite[Thm.\ 3.8.2]{ss-classification}.

A spectral model category is in particular a {\em simplicial} 
and {\em stable} model category.
Moreover, for $X$ a cofibrant and $Y$ a 
fibrant object of a spectral model category $\C$  there is a natural
isomorphism of graded abelian groups
$\pi^s_* \Hom_{\C}(X,Y) \iso {\HoC}(X,Y)_*$.
These facts are discussed in Lemma 3.5.2 of \cite{ss-classification}.

For of an object $P$ in a spectral model category, the function spectrum 
of $\End_{\C}(P)=\Hom_{\C}(P,P)$ is naturally a ring spectrum;
the multiplication is a special case of the composition product
\[  \circ \ : \ \Hom_{\C}(Y,Z) \ \sm \ \Hom_{\C}(X,Y)
\ \to \ \Hom_{\C}(X,Z) \ . \]
Also via the composition pairing, the function symmetric spectrum 
$\Hom_{\C}(P,X)$ becomes a right module 
over the symmetric ring spectrum  $\End_{\C}(P)$ for any object $X$.
In order for the endomorphism ring spectrum $\End_{\C}(P)$ to have the
correct homotopy type, the object $P$ should be both cofibrant and fibrant. 
In that case, the ring of homotopy groups $\pi_*\End_{\C}(P)$ 
is isomorphic to ${\HoC}(P,P)_{\ast}$, the ring of graded self maps of $P$ in 
the homotopy category of~$\C$.
Moreover, the homotopy type of the endomorphism ring spectrum then depends only
on the weak equivalence type of the object 
(see \cite[Cor.\ A.2.4]{ss-classification}).
Note that this is not completely obvious since taking endomorphisms is
{\em not} a functor.

If $P$ is a cofibrant object of a spectral model category,
then the functor
\[ \Hom_{\C}(P,-) \, : \, {\C} \ \to \ \Mod\End_{\C}(P) \]
is the right adjoint of a Quillen adjoint functor pair,
see \cite[3.9.3 (i)]{ss-classification}. 
The left adjoint is denoted 
\begin{equation}\label{eq-left derived of sm P} 
-\sm_{\End_{\C}(P)}P \, : \, \Mod\End_{\C}(P)\ \to \ \C \ . \end{equation}
For spectral model categories, the proof of 
Theorem \ref{thm-main-one} is now straightforward,
and very analogous to the proofs of Theorem \ref{thm-Gabriel/Michell} and 
Theorem \ref{main theorem DG version}; indeed, to obtain the
following theorem, one applies Proposition \ref{triangulated comparison prop}
to the total left derived functor
of the left Quillen functor \eqref{eq-left derived of sm P}.

\begin{theorem} \label{main theorem spectral version} Let $\C$  
be a spectral model category and  $P$ a cofibrant-fibrant object.
If $P$ is a  small generator for $\C$, then the  
adjoint functor pair $\Hom_{\C}(P,-)$ and $-\sm_{\End_{\C}(P)}P$ form a  
spectral Quillen equivalence.
\end{theorem}

The remaining step is worked out in Theorem 3.8.2 of 
\cite{ss-classification}, which proves that
every  simplicial, cofibrantly generated, proper stable model category 
is Quillen equivalent to a spectral model category,
namely the category $\spec(\C)$ of {\em symmetric spectra over $\C$}.
The proof is technical and we will not go into details here.
Theorem \ref{thm-main-one} follows by combining 
\cite[Theorem 3.8.2]{ss-classification} and 
Theorem \ref{main theorem spectral version} 
to obtain a diagram of model categories and  
Quillen equivalences (the left adjoints are on top)
\[\xymatrix@=20mm{ \C \quad \ar@<.4ex>^-{\Sigma^{\infty}}[r] & 
\quad  \spec(\C) \quad  \ar@<.4ex>^-{\text{Ev}_0}[l] 
\ar@<-.4ex>_-{\Hom_{\C}(P,-)}[r] & 
\quad  \Mod\End_{\C}(P) \ . \ar@<-.4ex>_-{- \sm_{\End_{\C}(P)}P}[l]} 
\]

\subsection{Morita context for ring spectra} \label{sec-Morita}

Now we come to `Morita theory for ring spectra', by which we mean the
question when two symmetric spectra have Quillen equivalent
module categories.
For ring spectra, there is a significant difference between
a Quillen equivalence of the module categories and an equivalence of the
homotopy categories. The former implies the latter, but not conversely.
The same kind of difference already exists for differential graded rings,
but it is not visible for ordinary rings
(see Example \ref{single generator examples} (5)).

We call a symmetric spectrum $X$ {\em flat} if the functor
$X\sm -$ preserves stable equivalences of symmetric spectra.
If $X$ is cofibrant, or more generally $S$-cofibrant 
in the sense of \cite[5.3.6]{hss}, then $X$ is flat, see~\cite[5.3.10]{hss}.
Every symmetric ring spectrum has a `flat resolution':
we may take a cofibrant approximation in the stable model structure
of symmetric ring spectra~\cite[5.4.3]{hss}; 
the underlying symmetric spectrum of the approximation is cofibrant, thus flat.

\begin{theorem} \label{thm-Morita} {\bf (Morita context)} 
The following are equivalent for two symmetric ring spectra $R$ and $S$.\\
\hspace*{0.5cm}{\em (1)} There exists a chain of spectral Quillen equivalences 
between the categories of $R$-modules and $S$-modules.\\
\hspace*{0.5cm}{\em (2)} There is a small, 
cofibrant and fibrant generator of the model category of $S$-modules
whose endomorphism ring spectrum is stably equivalent to $R$.\\

Both conditions are implied by the following condition.\\
\hspace*{0.5cm}{\em (3)} There exists an $R$-$S$-bimodule $M$ 
such that the derived smash product functor 
\[ - \sm^L_R M \ : \ \Ho(\Mod R) \ \to \ \Ho(\Mod S) \]
is an equivalence of categories.\\

If moreover $R$ or $S$ is flat as a symmetric spectrum, 
then all three conditions are equivalent.\\
\end{theorem}

Again there is a version of the Morita context \ref{thm-Morita}
relative to a commutative  symmetric ring spectrum $k$. 
In that case, $R$ and $S$ are $k$-algebras, condition (1) 
refers to $k$-linear spectral Quillen equivalences,
condition (2) requires a stable equivalence of $k$-algebras,
the bimodule $M$ in (3) has to be $k$-symmetric  
and in the addendum, one of $R$ or $S$ has to be flat as a $k$-module.

\begin{proof}[Proof of Theorem \ref{thm-Morita}]
{\bf (2)$\Longrightarrow$ (1):}
Modules over a symmetric ring spectrum form a spectral model
category; so this implication is a special case of
Theorem \ref{main theorem spectral version}, combined with the
fact that stably equivalent ring spectra have Quillen equivalent 
module categories.

{\bf (1) $\Longrightarrow$ (2):} 
To simplify things we suppose that there exists a single 
spectral Quillen equivalence
\[\xymatrix@=20mm{ \Mod R \quad \ar@<.4ex>^-{\Lambda}[r] & 
\quad  \Mod S \quad  \ar@<.4ex>^-{\Phi}[l] }\]
with  $\Lambda$ the left adjoint.
The general case of a chain of such Quillen equivalences 
is treated in~\cite[Thm.\ 4.1.2]{ss-classification}.
We choose a trivial cofibration $\iota:\Lambda(R) \to \Lambda(R)^{\text{f}}$
of $S$-modules such that $M:=\Lambda(R)^{\text{f}}$ is fibrant;
since $M$ is isomorphic in the homotopy category 
of $S$-modules to the image of the free $R$-module 
of rank one under the equivalence of homotopy categories,
$M$ is a small generator for the homotopy category of $S$-modules. 
It remains to show that the endomorphism
ring spectrum of $M$ is stably equivalent to $R$.

We define $\End_S(\iota)$, the endomorphism ring spectrum of the 
$S$-module map $\iota:\Lambda(R) \to M$, as the pullback in the diagram 
of symmetric spectra
\begin{equation}\label{eq-endo of map}
\xymatrix{ \End_S(\iota) \ar[r] \ar[d] & \End_S(M) \ar[d]^{\iota^*} \\
R \ar[r]_-{\iota_*} & \Hom_S(\Lambda(R),M) }\end{equation}
The right vertical map $\iota^*$ is obtained by applying $\Hom_S(-,M)$
to the acyclic cofibration $\iota$; since $M$ is stably fibrant,
$\iota^*$ is acyclic fibration of symmetric spectra.
Since $\Lambda$ and $\Phi$ form a Quillen equivalence, 
the adjoint $\hat \iota:R\to \Phi(M)$ of $\iota$ is a stable equivalence 
of $R$-modules. 
The lower horizontal map $\iota_*$ is the composite stable equivalence
\[ R \ \iso \ \Hom_R(R,R) \ \xrightarrow{\, \Hom_R(R,\hat\iota)\ } \ 
\Hom_R(R,\Phi(M)) \ \iso \ \Hom_S(\Lambda(R),M) \ . \] 
All maps in the pullback square \eqref{eq-endo of map}
are thus stable equivalences and the morphism connecting  $\End_S(\iota)$
to $R$ and  $\End_S(M)$ are homomorphisms of symmetric ring spectra
(whereas the lower right corner $\Hom_S(\Lambda(R),M)$ has no multiplication).
So $R$ is indeed stably equivalent, as a symmetric ring spectrum,
to the endomorphisms of $M$. 

{\bf (3)$\Longrightarrow$ (1):} If $M$ happens to be cofibrant as a right
$S$-module, then smashing with $M$ over $R$ is a left Quillen equivalence
from $R$-modules to $S$-modules. Since we did not assume that $M$
is cofibrant over $S$, we have to be content with a chain of two 
Quillen equivalences, which we get as follows.

Let $M$ be an $R$-$S$-bimodule as in condition (3). 
We choose a cofibrant approximation $\iota:R^{\text{c}}\to R$
in the stable model structure of symmetric ring spectra
and we view $M$ as an $R^{\text{c}}$-$S$-bimodule by restriction of scalars.
Then we choose a cofibrant approximation 
$M^{\text{c}}\to M$ as  $R^{\text{c}}$-$S$-bimodules.
Since the underlying symmetric spectrum of $R^{\text{c}}$ 
is cofibrant~\cite[4.1 (3)]{ss-monoidal},
$R^{\text{c}}\sm S$ is cofibrant as a right $S$-module, and thus
every cofibrant $R^{\text{c}}$-$S$-bimodule is cofibrant 
as a right $S$-module.
In particular, this holds for $M^{\text{c}}$, and so
we have a chain of two spectral Quillen pairs
\[\xymatrix@=20mm{ \Mod R \quad \ar@<-.4ex>_-{\iota^*}[r]  & 
\quad  \Mod R^{\text{c}} \quad   \ar@<-.4ex>_-{-\sm_{R^{\text{c}}}R}[l]
 \ar@<.4ex>^-{- \sm_{ R^{\text{c}}}M^{\text{c}}}[r] & 
\quad  \Mod S \ . \ar@<.4ex>^-{\Hom_S(M^{\text{c}},-)}[l]  } 
\]
The left pair is a Quillen equivalence since the approximation map 
$\iota:R^{\text{c}}\to R$ is a stable equivalence. 
For every cofibrant $R^{\text{c}}$-module $X$, the map
\[ X\sm_{R^{\text{c}}}M^{\text{c}} \ \to \
X\sm_{R^{\text{c}}}M \iso (X\sm_{R^{\text{c}}}R)\sm_R M \]
is a stable equivalence.
This means that the diagram of homotopy categories and derived functors
\[\xymatrix{ & \Ho(\Mod R^{\text{c}}) \ar[dl]_{-\sm^L_{R^{\text{c}}}R}
\ar[dr]^{-\sm^L_{R^{\text{c}}}M^{\text{c}}} \\
\Ho(\Mod R) \ar[rr]_{-\sm^L_{R}M} && \Ho(\Mod S) }  \] 
commutes up to natural isomorphism.
Thus the right Quillen pair above induces 
an equivalence of homotopy categories, so it is a Quillen equivalence.

{\bf (2)$\Longrightarrow$(3), assuming that $R$ or $S$ is flat.}
Let $T$ be a cofibrant and fibrant small generator of $\Ho(\Mod S)$ 
such that $R$ is stably equivalent to the endomorphism ring spectrum of $T$.
We choose a cofibrant approximation $R^{\text{c}}\xrightarrow{\simeq} R$ in the
model category of symmetric ring spectra. Since $T$ is cofibrant and fibrant,
its endomorphism ring spectrum is fibrant. So any isomorphism
between $R$ and $\End_S(T)$ in the homotopy category of symmetric ring
spectra can be represented by a chain of two stable equivalences
\[ R \ \xleftarrow{\ \simeq\ } \ R^{\text{c}} \ \xrightarrow{\ \simeq\ } \
\End_S(T) \ .  \]
The module $T$ is naturally an $\End_S(T)$-$S$-bimodule,
and we restrict the left action to $R^{\text{c}}$ and view $T$ as an
$R^{\text{c}}$-$S$-bimodule.
We choose  a cofibrant replacement $T^{\text{c}}\stackrel{\sim}{\to}T$ 
as an $R^{\text{c}}$-$S$-bimodule. Then we set 
\[ M \ = \ R\sm_{R^{\text{c}}} T^{\text{c}} \ , \] 
an $R$-$S$-bimodule.
We have no reason to suppose that $M$ is cofibrant as a right $S$-module,
so we cannot assume that the functor $-\sm_RM:\Mod R\to \Mod S$ is
a left Quillen functor. Nevertheless, smashing with $M$ over $R$
takes stable equivalences between cofibrant $R$-modules
to stable equivalences, so it has a total left derived functor
\[ -\sm_R^L M \ : \ \Ho(\Mod R) \ \to \ \Ho(\Mod S) \ ; \] 
we claim that this functor is an equivalence.

Since $R^{\text{c}}$ is cofibrant as a symmetric ring spectrum, it is also
cofibrant as a symmetric spectrum \cite[4.1 (3)]{ss-monoidal}, 
so $R^{\text{c}}\sm S^{op}$ models the derived smash product of $R$ and $S$. 
If one of $R$ or $S$ are flat, then $R\sm S^{op}$ 
also models the derived smash  product, so that the map
\[ R^{\text{c}}\sm S^{op} \ \to \ R\sm S^{op} \] 
is a stable equivalence of symmetric ring spectra. 
Since $T^{\text{c}}$ is cofibrant as an $R^{\text{c}}\sm S^{op}$-module, 
the induced map
\begin{equation}\label{eq-an equivalence}  
T^{\text{c}} \ = \  (R^{\text{c}}\sm S^{op})_{R^{\text{c}}\sm S^{op}} T^{\text{c}} \ \to \ 
 (R\sm S^{op})_{R^{\text{c}}\sm S^{op}} T^{\text{c}} \ \iso \
R\sm_{R^{\text{c}}} T^{\text{c}} \ = \ M \end{equation}
is a stable equivalence.
We smash the  stable equivalence \eqref{eq-an equivalence}  
from the left with an $R^{\text{c}}$-module $X$ to get a natural map 
of $S$-modules
\begin{equation}\label{eq:two maps}  
 X\sm_{R^{\text{c}}} T^{\text{c}} \
\xrightarrow{\ }  \   X\sm_{R^{\text{c}}}M \iso  
(X\sm_{R^{\text{c}}}R)\sm_R M  \ . \end{equation} 
If $X$ is cofibrant as an $R^{\text{c}}$-module, 
then $X\sm_{R^{\text{c}}}-$ takes
stable equivalences of left $R^{\text{c}}$-modules to stable equivalences, 
so in this case, the map \eqref{eq:two maps} is a stable equivalence.
Thus the diagram of triangulated categories and derived functors
\begin{equation}\label{eq-claimed triangle}
\xymatrix{ 
& \Ho(\Mod R^{\text{c}}) \ar[dr]^-{-\sm^L_{R^{\text{c}}} T^{\text{c}}}  
\ar[dl]_{-\sm^L_{R^{\text{c}}} R}\\
\Ho(\Mod R) \ar[rr]_-{ -\sm_R^L M} && \Ho(\Mod S) } \end{equation}
commutes up to natural isomorphism.

The left diagonal functor in the diagram \eqref {eq-claimed triangle}
is derived from extensions of scalars
along stable equivalences of ring spectra; such extension of scalars
is a left Quillen equivalences, so the derived functor 
$-\sm^L_{R^{\text{c}}} R$ is an exact equivalence of triangulated categories.
We argued in the previous implication that any cofibrant  
$R^{\text{c}}\sm S^{op}$-module such as $T^{\text{c}}$ 
has an underlying cofibrant right $S$-module.
So smashing with $T^{\text{c}}$ over $R^{\text{c}}$ is a left Quillen functor. 
Since $T^{\text{c}}$ is isomorphic to $T$ in the homotopy category of
$S$-modules, $T^{\text{c}}$ is a  small generator of $\Ho(\Mod S)$.
So the right diagonal derived functor
in \eqref {eq-claimed triangle}is an exact equivalence
by Proposition \ref{triangulated comparison prop}, applied to the
free $R^{\text{c}}$-module of rank 1.
So we conclude that the lower horizontal functor 
in the diagram \eqref {eq-claimed triangle}
is also an exact equivalence of triangulated categories.
This establishes condition (3).
\end{proof}

\subsection{Examples}\smallskip
\label{single generator examples} 

(1) {\bf Matrix ring spectra.} 
As for classical rings (compare Example \ref{ex-classical matrix}),
matrix ring spectra give rise to the simplest kind of Morita equivalence.
Indeed over any a ring spectrum $R$, the `free module of rank $n$', i.e.,
the wedge of $n$ copies of $R$, is a small generator for the
homotopy category of $R$-modules. The endomorphism ring spectrum of 
a (stably fibrant replacement) of $R\times n$ is the $n\times n$ matrix
ring spectrum as we defined it in Example \ref{ex-matrix ring spectrum}.
So $R$ and 
\[ M_n(R) \ = \ \End_{R}((R \times n)^{\text{f}}) \] 
are Morita equivalent as ring spectra.

(2) {\bf Upper triangular matrices.}
In Example \ref{ex-tilting module A3 quivers} we saw that the
upper triangular $3\times 3$ matrices over a field are
derived equivalent, but not Morita equivalent, to its
sub-algebra of matrices of the form
\[ \left\lbrace 
\begin{pmatrix}   x_{11} & x_{12} & x_{13} \\
                  0  & x_{22} & 0 \\
                  0  &   0  & x_{33} \end{pmatrix} \ | \
x_{ij} \in k \right\rbrace  \ . \]
We cannot directly define the algebra of  $3\times 3$ matrices over a
ring spectrum; the problem is that the usual basis of elementary matrices
is closed under multiplication, but the unit matrix
is a {\em sum} basis elements, not a single basis element.

The algebras of  Example \ref{ex-tilting module A3 quivers} 
can also be obtained as the path algebras of the two quivers
\[ Q=\left\lbrace 
\xymatrix{ 1 \ar[r]^{\alpha} & 2 \ar[r]^{\beta} & 3 } 
\right\rbrace \ \text{\quad respectively\quad}
\left\lbrace \xymatrix{ 1 & 2 \ar[l]  \ar[r] & 3 } 
\right\rbrace \ . \]
With this in mind we can now run  Example \ref{ex-tilting module A3 quivers} 
with the ground field replaced by a commutative symmetric ring spectrum $R$.

We construct the `path algebra' indirectly via 
`representations of the quiver $Q$ over $R$'. A representation of $Q$ over $R$
is a collection
\[ M \ = \ \left\lbrace \, M_1 \xrightarrow{\ \alpha\ } M_2 
\xrightarrow{\ \beta\ } M_3 \, \right\rbrace \]
of three $R$-modules and two homomorphisms.
A morphisms $f:M\to N$ of representations consists 
of $R$-homomorphisms $f_i:M_i\to N_i$ for $i=1,2,3$ satisfying
$\alpha f_1=f_2\alpha$ and $\beta f_2=f_3\beta$.
There is a stable model structure on the category of such representations
in which a morphism  $f:M\to N$ is a stable equivalence or fibration 
if and only if each $R$-homomorphism $f_i:M_i\to N_i$ is a stable equivalence 
or fibration for all $i=1,2,3$. 
Moreover, $f$ is a cofibration if and only if the morphisms
\[ f_1:M_1\to N_1 \, , \ f_2\cup\alpha : M_2\cup_{M_1}N_1 \to N_2
\quad \text{and} \quad  f_3\cup\beta : M_3\cup_{M_2}N_2 \to N_3 \]
are cofibrations of $R$-modules. 

We consider the `free' or `projective' representations $P^i$
given by $P^1=\{R\stackrel{=}{\to}R\stackrel{=}{\to}R\}$, 
$P^2=\{\ast\to R\stackrel{=}{\to}R\}$ respectively $P^3=\{\ast\to\ast\to R\}$.
The representation $P^i$ represents the evaluation functor
functor, i.e., we have a natural isomorphism
\[ \Hom_{Q\text{-rep}}(P^i,M) \ \iso \ M_i \ .  \]
This is implies that the wedge of these three representations
is a small generator of the homotopy category of $Q$-representations.
The three projective representations are cofibrant;
we let $M_3^{\triangle}(R)$ denote the endomorphism ring spectrum
of a stably fibrant approximation of their wedge,
\[  M_3^{\triangle}(R) \ = \  
\End_{Q\text{-rep}}((P^1\vee P^2\vee P^3)^{\text{f}}) \ . \] 
This symmetric ring spectrum deserves to be called the 
`upper triangular $3\times 3$ matrices' over $R$.

Now we find a different small generator for the model category
of $Q$-representations which is the analog of the tilting module
in Example \ref{ex-tilting module A3 quivers}.
We note that the inclusion $P^3\to P^2$ is not a cofibration, and the quotient
$P^2/P^3=\{\ast\to R\to\ast\}$ is not cofibrant. A cofibrant approximation
of this quotient is given by the representation
$S^2=\{\ast\to R\to CR\}$ in which the second map is the inclusion
of $R$ into its cone. We form the representation
\[ T \ = \ P^1 \, \vee \, P^2 \, \vee \, S^2 \ . \] 
Then $S^2$ represents the functor which sends
a $Q$-representation $M$ to the homotopy fiber of the morphism
$\beta:M_2\to M_3$, and thus
\[  \Hom_{Q\text{-rep}}(T,M) \ \iso \ M_1 \times M_2 \times
\text{hofibre}(\beta:M_2\to M_3) \ . \]  
This implies that $T$ is also a small generator
for the homotopy category of $Q$-representation over $R$.

We conclude that the upper triangular matrix algebra
$M_3^{\triangle}(R)$ is Quillen equivalent to the endomorphism
ring spectrum of a stably fibrant approximation of 
the representation $T$. Just as the endomorphisms of the
generator $P^1\vee P^2\vee P^3$ should be thought of
as upper triangular matrices, the endomorphisms of the
generator $P^1\vee P^2\vee S^2$ are analogous to a certain algebra of
$3\times 3$ matrices over $R$, namely the ones of the form
\[ \begin{pmatrix}     R  &   R  &  R   \\
                     \ast &   R  & \ast \\
                     \ast & \ast &  R \end{pmatrix}   \ . \]
Another example is obtained as follows. We consider
the representation $S^1=\{R\to CR\stackrel{=}{\to} CR\}$
which is a cofibrant replacement of the representation $P^1/P^2$
and which represents the homotopy fiber of the morphism
$\alpha:M_1\to M_2$. Then 
\[ T' \ = \  S^1 \, \vee \, S^2 \, \vee \, P^3 \] 
is another small generator for the homotopy category 
of $Q$-representation over $R$. So $M_3^{\triangle}(R)$ 
is also Quillen equivalent to the derived endomorphism
ring spectrum of $T'$, which is an algebra of
$3\times 3$ matrices of the form
\[ \begin{pmatrix}     R  &  \Omega R  &    \ast   \\
                     \ast &      R     &  \Omega R \\
                     \ast &    \ast    &      R \end{pmatrix}   \ , \]
where each symbol `$\ast$' designates an entry
in a stably contractible spectrum.

(3) {\bf Uniqueness results for stable homotopy theory.}
Theorem \ref{thm-main-one} characterizes module categories over
ring spectra among stable model categories.
This also yields a characterization of the model category of spectra: 
a stable model category is Quillen equivalent
to the category of symmetric spectra if and only if it has 
a small generator $P$ for which the unit map of ring spectra
$\mS\to\End(P)$ is a stable equivalence.
The technical conditions of being a simplicial, cofibrantly generated
and proper can be eliminated as in  \cite{ss-uniqueness} 
with the use of {\em framings}~\cite[Chpt.~5]{hovey-book}.
The paper \cite{ss-uniqueness} also gives other necessary and 
sufficient conditions for when a stable model category is 
Quillen equivalent to spectra -- some of them in terms of the
homotopy category and the natural action of the stable homotopy
groups of spheres. 
In \cite{sch-2}, this result is extended to a uniqueness theorem
showing that the 2-local stable homotopy category has only one
underlying model category up to Quillen equivalence;
the odd-primary version is work in progress.
In another direction, the uniqueness result is extended to include the
monoidal structure in \cite{shipley-mon}.

(4) {\bf Chain complexes and Eilenberg-Mac Lane spectra.} 
For a ring $R$, the category of chain complexes of $R$-modules 
(under quasi-isomorphisms) is Quillen  equivalent to the category 
of modules over the Eilenberg-Mac Lane ring spectrum $HR$.
An explicit chain of two Quillen equivalences can be found in 
Theorem B.1.11 of \cite{ss-classification};
I don't know if it is possible to compare the two categories 
by a single Quillen equivalence.

This result can also be viewed as an instance of Theorem \ref{thm-main-one}: 
the free $R$-module of rank one, considered as a complex concentrated in 
dimension zero, is a small generator for the unbounded derived category 
of $R$.
Since the homotopy groups of its endomorphism ring spectrum 
(as an object of the model category of chain complexes) are concentrated
in dimension zero, the endomorphism ring spectrum is stably equivalent to the 
Eilenberg-Mac Lane ring spectrum for $R$ (Proposition B.2.1 of 
\cite{ss-classification}).

(5) {\bf A generalized tilting theorem.}
We interpret and generalize the tilting theory 
from the perspective of stable model categories. 
A {\em tilting object} in a stable model category $\C$
as a small generator $T$ such that the graded homomorphism group
$[T,T]_*^{\HoC}$ in the homotopy category is concentrated in dimension zero. 
The following `generalized tilting theorem' 
of \cite[Thm. 5.1.1]{ss-classification} then
shows that the existence of a tilting object is necessary and sufficient 
for a stable model category to be Quillen equivalent 
or derived equivalent to the category of unbounded chain complexes over a ring.
\vspace*{.2cm}

{\bf Generalized tilting theorem.}
{\em Let $\C$ be a stable model category which is simplicial,
cofibrantly generated and proper, and let $R$ be a ring.
 Then the following conditions are equivalent:\\
\hspace*{0.5cm}{\em (1)} There is a chain of Quillen equivalences between $\C$ 
and the model category of chain complexes of $R$-modules.\\
\hspace*{0.5cm}{\em (2)} The homotopy category of $\C$ is 
triangulated equivalent to the derived category $\D(R)$.\\
\hspace*{0.5cm}{\em (3)} The model category $\C$ has a tilting object
whose endomorphism ring in $\Ho(\C)$ is isomorphic to $R$.
\vspace*{.2cm}}

In the derived category of a ring, a tilting object is the same as a
tilting complex, and the result reduces to 
Rickard's tilting theorem \ref{thm-tilting}. 

The generalized tilting situation enjoys one very special feature.
In general, the notion of Quillen equivalence 
is considerably stronger than triangulated equivalence of homotopy categories;
two examples are given in \cite[2.1 and 2.2]{sch-2}.
Hence it is somewhat remarkable that for complexes of modules 
over rings, the two notions are in fact equivalent. 
In general the homotopy category determines the homotopy groups of the 
endomorphism ring spectrum, but not its homotopy type. 
The real reason behind the equivalences of conditions (1) and (2) 
above is the fact that in contrast to arbitrary ring spectra,
Eilenberg-Mac Lane spectra are determined up to stable equivalence
by their homotopy groups. We explain this in more detail in Section 5 of
\cite{ss-classification}.

(6) {\bf Frobenius rings.}
As in Example \ref{ex-stable model category examples} (2) we consider  
a Frobenius ring and assume that the stable module category has a  
small generator. Then we are in the situation of Theorem  
\ref{thm-main-one}; however this example is completely algebraic,  
and there is no need to consider ring spectra to identify the stable  
module category as the derived category of a suitable `ring'. 
In fact Keller shows \cite[4.3]{keller-derivingDG} that in such a  
situation there exists a differential graded algebra and  
an equivalence between the stable module category and the unbounded  
derived category of the differential graded algebra. 

(7) {\bf Smashing Bousfield localizations.} 
Let $E$ be a spectrum and consider the $E$-local model  
category structure on some model category of spectra 
(see e.g.\ \cite[VIII 1.1]{ekmm}).
This is another stable model category in which the localization of  
the sphere spectrum $L_E \mS$ is a generator. This localized sphere  
is small if the localization is {\em smashing}, i.e., if a certain  
natural map $X\sm L_E \mS \to L_E X$ is a stable equivalence for  
all $X$. So for a smashing localization the $E$-local model category  
of spectra is Quillen equivalent to modules over the ring spectrum  
$L_E \mS$ (which is the endomorphism ring spectrum of the localized  
sphere in the localized model structure).

(8) {\bf Finite localization.}
\label{rem-fin-loc} 
Suppose $P$ is a small object of a triangulated category $\T$ 
with infinite coproducts. Then there always exists an idempotent 
localization functor $L_P$ on $\T$ whose acyclics are precisely 
the objects of the localizing subcategory generated by $P$ 
(compare \cite{miller-finite} or the proofs of
\cite[Lemma 2.2.1]{ss-classification} or 
\cite[Prop.\ 2.3.17]{hps}). These localizations are often referred to 
as {\em finite Bousfield localizations} away from $P$. 

This type of localization has a refinement to the model category level. 
Suppose $\C$ is a stable model category and $P$ 
a small object, and let $L_P$ denote the associated
localization functor on the homotopy category of $\C$. 
By \ref{main theorem spectral version},
or rather the refined version \cite[Thm.\ 3.9.3 (ii)]{ss-classification}, 
the acyclics for $L_P$ are equivalent 
to the homotopy category of $\End_{\C}(P)$-modules, 
the equivalence arising from a Quillen adjoint functor pair. 
Furthermore the counit of the derived adjunction
\[ \Hom_{\C}(P,X) \, \sm_{\End_{\C}(P)}^L \, P \ \to \ X \]
is the acyclicization map and its cofiber is a model for the 
localization $L_P X$.

(9) {\bf $K(n)$-local spectra.} Even if a  
Bousfield localization is not smashing, Theorem \ref{thm-main-one}  
might be applicable. As an example we consider Bousfield  
localization with respect to the $n$-th Morava K-theory $K(n)$ at a  
fixed prime.
The localization of the sphere is still a generator,
but for $n>0$
it is not small in the local category, see~\cite[3.5.2]{hps}.
However the localization of any finite type $n$ spectrum $F$ is a  
small generator for the $K(n)$-local  
category~\cite[7.3]{hovey-strickland-kn}. Hence the $K(n)$-local  
model category is Quillen
equivalent to modules over the endomorphism ring spectrum of
$L_{K(n)}F$.

\medskip

I would like to conclude with a few words about {\bf invariants} 
of ring spectra which are preserved under Quillen equivalence
(but {\em not} in general under equivalences of 
triangulated homotopy categories).
Such invariants include the algebraic $K$-theory~\cite{dugger-shipley}, 
topological Hochschild homology and topological cyclic homology.

In the classical framework, the center of a ring is invariant
under Morita and derived equivalence. As a general philosophy 
for spectral algebra, definitions which use elements 
are not well-suited for generalization to ring spectra.
So how do we define the `center' of a ring spectrum, such that it only depends,
up to stable equivalence, on the Quillen equivalence class 
of the module category ? 
The center of an ordinary ring $R$ is isomorphic to the
endomorphism ring of $R$, considered as a bimodule over itself, via
\[ \End_{R\tensor R^{\text{op}}}(R) \ \xrightarrow{\quad } \ Z(R) \ , 
\quad f \longmapsto f(1) \]
So we define the {\em center} of a ring spectrum $R$ 
as the endomorphism ring spectrum
of a cofibrant-fibrant replacement of $R$, considered as a bimodule
over itself,
\begin{equation}\label{brave new center} 
Z(R) \ = \ \End_{R\sm R^{\text{op}}}(R^{\text{cf}},R^{\text{cf}})\ .  
\end{equation}
In this definition, $R$ should be flat as a symmetric spectrum,
in order for the smash product $R\sm R^{\text{op}}$ to have 
the `correct homotopy type'.

This kind of center of a ring spectrum is homotopy commutative; 
slightly more is true:
the center~\eqref{brave new center} is often called the 
{\em topological Hochschild cohomology spectrum}
of the ring spectrum $R$, and its multiplication extends to
an action of an operad weakly equivalent 
to the operad of little discs, see \cite[Thm.\ 3]{mcclure-smith:deligne}.
However, the above center is  usually {\em not} stably equivalent
to a {\em commutative} symmetric ring spectrum (or what is the same,
$E_{\infty}$-homotopy commutative);
the Gerstenhaber operations on the homotopy of $Z(R)$ are obstruction to
higher order commutativity. So is it just a coincidence that the classical
center of an ordinary ring is commutative ? Or is there some
`higher', $E_{\infty}$-commutative, center of a ring spectrum,
yet to be discovered ?


\begin{thebibliography}{EKMM}

\bibitem[AF92]{anderson-fuller}
F.~W.~Anderson and K.~R.~Fuller, {\em Rings and categories of modules 
(second edition)},
Graduate Texts in Mathematics, {\bf 13}, Springer-Verlag, New York, 1992,
viii+376 pp.

\bibitem[AGH]{Morita obituary}
A.~V.~Arhangel'skii, K.~R.~Goodearl, and B.~Huisgen-Zimmermann,
{\em Kiiti Morita (1915--1995)},
Notices Amer. Math. Soc. 44 (1997), no. 6, 680--684.

\bibitem[Ba]{bass-Ktheory}
H.~Bass, {\em Algebraic $K$-theory.} W.~A.~Benjamin, Inc., 
New York-Amsterdam 1968 xx+762 pp. 

\bibitem[BL]{BL}
A.~Baker, A.~Lazarev, {\em Topological Hochschild cohomology and 
generalized Morita equivalence}, Preprint (2002).
{\verb|http://arXiv.org/abs/math.AT/0209003|}

\bibitem[BBD]{BBD}
A.~A.~Beilinson, J.~Bernstein, P.~Deligne,
{\em Faisceaux pervers.} Analysis and topology on singular spaces, I 
(Luminy, 1981), 5--171, Ast\'erisque {\bf 100}, 
Soc.\ Math.\ France, Paris, 1982. 

\bibitem[Bei78]{beilinson}
A.~A.~Beilinson, {\em Coherent sheaves on $ P\sp{n}$ and 
problems in linear algebra},
Functional Anal. Appl. 12 (1978), no. 3, 214--216 (1979).

\bibitem[Bek00]{beke-sheafy}
T.~Beke, {\em Sheafifiable homotopy model categories.}
Math. Proc. Cambridge Philos. Soc. {\bf 129} (2000), 447--475.

\bibitem[Bou75]{bousfield-localization-spaces}
A.~K.~Bousfield, {\em The localization of spaces with respect to homology}, 
Topology {\bf 14} (1975), 133--150.

\bibitem[BF78]{BF}
A.~K.~Bousfield and E.~M.~Friedlander, {\em Homotopy theory of 
{$\Gamma$}-spaces, spectra, and bisimplicial sets}, Geometric applications of 
homotopy theory (Proc. Conf., Evanston, Ill., 1977), II (M.~G. Barratt and 
M.~E. Mahowald, eds.), Lecture Notes in Math., {\bf 658}, Springer, Berlin, 
1978, pp.~80--130.

\bibitem[Bro94]{broue-block equivalences}
M.~Brou{\'e}, {\em Equivalences of blocks of group algebras.}
Finite-dimensional algebras and related topics (Ottawa, ON, 1992), 1--26, 
NATO Adv. Sci. Inst. Ser. C Math. Phys. Sci., 424, 
Kluwer Acad. Publ., Dordrecht, 1994. 

\bibitem[CH]{Christensen-Hovey-relative}
D.~Christensen, M.~Hovey, 
{\em Quillen model structures for relative homological algebra},
Math.\ Proc.\ Cambridge Philos.\ Soc. {\bf 133} (2002), 261--293.

\bibitem[SGA 4${1\over 2}$]{SGA 4 1/2}
P.~Deligne, {\em Cohomologie \'etale}. S\'eminaire de G\'eom\'etrie 
Alg\'ebrique du Bois-Marie SGA 4${1\over 2}$. 
Avec la collaboration de J.~F.~Boutot, A.~Grothendieck, L.~Illusie 
et J.~L.~Verdier. Lecture Notes in Mathematics, Vol. 569.
Springer-Verlag, Berlin-New York, 1977. iv+312pp. 

\bibitem[Dug]{dugger-simplicial}
D.~Dugger, {\em Replacing model categories with simplicial ones}, 
Trans.\ Amer.\ Math.\ Soc.\ {\bf 353} (2001), 5003--5027.

\bibitem[DuSh]{dugger-shipley}
D.~Dugger and B.~Shipley, {\em $K$-theory and derived equivalences},
Preprint (2002)\\ 
{\verb|http://arXiv.org/abs/math.KT/0209084|}

\bibitem[D{\O}R]{DOR}
B.~Dundas, P.~A.~{\O}stv{\ae}r and O.~R{\"o}ndigs, 
{\em Enriched functors and motivic stable homotopy theory.}
Preprint (2002) \verb!http://www.math.ntnu.no/~dundas/!

\bibitem[DwSp95]{DS}
W.~G.~Dwyer and J.~Spalinski, {\em Homotopy theories and model categories}, 
Handbook of algebraic topology (Amsterdam), North-Holland, 
Amsterdam, 1995, pp.~73--126.

\bibitem[DG02]{DG}
W.~G.~Dwyer, J.~Greenlees, {\em Complete modules and torsion modules},
Amer. J. Math. 124 (2002), 199--220.

\bibitem[DGI]{DGI}
W.~G.~Dwyer, J.~Greenlees, S.~Iyengar. 
{\em Duality in algebra and topology}, Preprint (2002)\\
\verb!http://www.shef.ac.uk/~pm1jg/!

\bibitem[EKMM]{ekmm}
A.~D.~Elmendorf, I.~Kriz, M.~A.~Mandell, and J.~P.~May, 
{\em Rings, modules, and algebras in stable homotopy theory. 
{W}ith an appendix by M.~Cole},
Mathematical Surveys and Monographs, {\bf 47}, American Mathematical Society,
Providence, RI, 1997, xii+249 pp.

\bibitem[GM]{Gelfand-Manin}
S.~I.~Gelfand, Y.~I.~ Manin, {\em Methods of homological algebra}. 
Translated from the 1988 Russian original. 
Springer-Verlag, Berlin, 1996. xviii+372 pp.

\bibitem[Gr77]{grothendieck-groupes des classes}
A.~Grothendieck, {\em Groupes des classes des cat{\'e}gories ab\'eliennes
et triangul{\'e}es, complexes parfaits.} Expos\'e VIII,
S\'eminaire de G\'eom\'etrie Alg\'ebrique du Bois-Marie SGA 5,
Lecture Notes in Math., {\bf 589}, Springer, Berlin, 
1977, pp.~351--371.
 
\bibitem[Hir03]{hirschhorn-book}
P.~S.~Hirschhorn, {\em Model Categories and Their Localizations},
Mathematical Surveys and Monographs {\bf 99},
American Mathematical Society, 2003, 457 pp.

\bibitem[Hov99]{hovey-book}
M.~Hovey, {\em Model categories}, Mathematical Surveys and Monographs,
{\bf 63}, American Mathematical Society,
Providence, RI, 1999, xii+209 pp.

\bibitem[Hov01a]{hovey-sheaves}
M.~Hovey, {\em Model category structures on chain complexes of sheaves}, 
Trans.\ Amer.\ Math.\ Soc.\ {\bf 353} (2001), 2441--2457.

\bibitem[Hov01b]{hovey-stabilization}
M.~Hovey, {\em Spectra and symmetric spectra in general model categories}, 
J.\ Pure Appl.\ Algebra {\bf 165} (2001), 63--127.

\bibitem[HPS97]{hps}
M.~Hovey, J.~H.~Palmieri, and N.~P.~Strickland, {\em Axiomatic stable
homotopy theory}, Mem. Amer. Math. Soc. {\bf 128} (1997), no.~610.

\bibitem[HSS]{hss}
M.~Hovey, B.~Shipley, and J.~Smith, {\em Symmetric spectra},
J.\ Amer.\ Math.\ Soc.\ {\bf 13} (2000), 149--208.

\bibitem[HS99]{hovey-strickland-kn}
M.~Hovey and N.~P.~Strickland, {\em Morava ${K}$-theories and localisation},
Mem. Amer. Math. Soc. {\bf 139} (1999)

\bibitem[Hu03]{hu-motivic S-modules}
P.~Hu, {\em $S$-modules in the category of schemes.}
Mem. Amer. Math. Soc. {\bf 161} (2003), no. 767, viii+125 pp.

\bibitem[Jar87]{jardine-stable_presheaves}
J.~F.~Jardine, {\em Stable homotopy theory of simplicial presheaves},
Canad. J. Math. {\bf 39} (1987), 733--747.

\bibitem[Jar97]{jardine-etale}
J.~F.~Jardine,  {\em Generalized \'etale cohomology theories},
Progress in Mathematics {\bf 146}, Birkh\"auser Verlag, Basel, 1997, x+317 pp.

\bibitem[Jar00a]{jardine-symmetric_presheaves}
J.~F.~Jardine, {\em Presheaves of symmetric spectra}, 
J.\ Pure Appl.\ Algebra {\bf 150} (2000), 137--154.

\bibitem[Jar00b]{jardine-A^1stable}
J.~F.~Jardine,  {\em  Motivic symmetric spectra},
Doc.\ Math.\ {\bf 5} (2000), 445--553.

\bibitem[Kel94a]{keller-derivingDG}
B.~Keller, {\em Deriving {D}{G} categories},
Ann. Sci. \'Ecole Norm. Sup. (4) {\bf 27} (1994), 63--102.

\bibitem[Kel94b]{keller-smashing}
B.~Keller, {\em A remark on the generalized smashing conjecture}, 
Manuscripta Math. {\bf 84} (1994), 193--198.

\bibitem[Kel96]{keller-cyclic invariance}
B.~Keller, {\em Invariance of cyclic homology under derived equivalence.} 
Representation theory of algebras (Cocoyoc, 1994), 353--361, 
CMS Conf. Proc., 18, Amer. Math. Soc., Providence, RI, 1996. 

\bibitem[Kel98]{keller-cyclic invariance and localization}
B.~Keller, {\em Invariance and localization for cyclic homology 
of DG algebras.} J. Pure Appl. Algebra {\bf 123} (1998), 223--273.

\bibitem[Kel99]{keller-cyclic of exact}
B.~Keller, {\em On the cyclic homology of exact categories}, 
J. Pure Appl. Algebra {\bf 136} (1999), 1--56.

\bibitem[Kel03]{keller-HochschildPicard}
B.~Keller, {\em Hochschild cohomology and derived Picard groups},
Preprint (2003)\\ 
{\verb|http://www.math.jussieu.fr/~keller/|}

\bibitem[KZ]{koenig-zimmermann}
S.~K\"onig, A.~Zimmermann, {\em Derived equivalences for group rings}.
With contributions by B.~Keller, M.~Linckelmann, J.~Rickard and R.~Rouquier. 
Lecture Notes in Mathematics, 1685.
Springer-Verlag, Berlin, 1998. x+246 pp.

\bibitem[KM]{kriz-may}
I.~Kriz and J.~P.~ May, {\em Operads, algebras, modules and motives}.
Ast{\'e}risque No.\ {\bf 233}, (1995), iv+145pp.

\bibitem[Lam]{Lam}
T.~Y.~Lam, {\em Lectures on modules and rings}. 
Graduate Texts in Mathematics, 189. Springer-Verlag, New York, 1999. 
xxiv+557 pp. 

\bibitem[LMS86]{lms}
L.~G.~Lewis, Jr., J.~P.~May, and M.~Steinberger, 
{\em Equivariant stable homotopy theory}, Lecture Notes in Mathematics, 
{\bf 1213}, Springer-Verlag, 1986.

\bibitem[Lyd98]{lydakis-simplicial}
M.~Lydakis, {\em Simplicial functors and stable homotopy theory},
Preprint (1998)\\
{\verb|http://hopf.math.purdue.edu/|}

\bibitem[McCS]{mcclure-smith:deligne}
J.~McClure, J.~Smith, 
{\em A solution of Deligne's Hochschild cohomology conjecture},
Recent progress in homotopy theory (Baltimore, MD, 2000), 153--193, 
Contemp. Math. {\bf 293}, Amer. Math. Soc., Providence, RI, 2002. 

\bibitem[McL]{MacLane-working}
S.~Mac Lane, {\em Categories for the working mathematician}. 
Graduate Texts in Mathematics, Vol.~5. Springer-Verlag, New York-Berlin, 1971.
 ix+262 pp. 

\bibitem[MM02]{mm}
M.~A.~Mandell and J.~P.~May, {\em Equivariant orthogonal spectra and
$S$-modules}, Mem. Amer. Math. Soc. {\bf 159} (2002), x+108 pp.

\bibitem[MMSS]{mmss}
M.~A.~Mandell, J.~P.~May, S.~Schwede and B.~Shipley,
{\em Model categories of diagram spectra}, 
Proc.\ London Math.\ Soc. {\bf 82} (2001), 441-512.

\bibitem[Mil92]{miller-finite}
H.~Miller, {\em Finite localizations}, 
Bol. Soc. Mat. Mexicana (2) {\bf 37} (1992),
(Papers in honor of Jos\'e Adem), 383--389.

\bibitem[MV]{morel-voevodsky}
F.~Morel and V.~Voevodsky, {\em {${\mathbb A}^1$}-homotopy theory of schemes},
Inst.~Hautes \'Etudes Sci.~Publ.~Math. {\bf 90} (2001), 45--143.
\bibitem[Mo58]{morita}
K.~Morita, {\em Duality for modules and its applications 
to the theory of rings with minimum condition.} 
Sci.\ Rep.\ Tokyo Kyoiku Daigaku Sect.\ A 6 (1958), 83--142. 

\bibitem[Nee92]{neeman-TTYBR}
A.~Neeman, {\em The connection between the $K$-theory localization theorem 
of Thomason, Trobaugh and Yao and the smashing subcategories 
of Bousfield and Ravenel.} 
Ann. Sci. \'Ecole Norm. Sup. (4) {\bf 25} (1992), 547--566.

\bibitem[Nee01]{neeman-triangle book}
A.~Neeman, {\em Triangulated categories}, 
Annals of Mathematics Studies, 148. 
Princeton University Press, Princeton, NJ, 2001. viii+449 pp.

\bibitem[Qui67]{Q}
D.~G. Quillen, {\em Homotopical algebra}, Lecture Notes in Mathematics,
{\bf 43}, Springer-Verlag, 1967.

\bibitem[Ren69]{rentschler}
R.~Rentschler, {\em Sur les modules $M$ tels que ${\rm Hom}(M,\,-)$ 
commute avec les sommes directes},
C.~R.~Acad. Sci. Paris S{\'e}r. A-B {\bf 268} (1969), A930--A933.

\bibitem[RSS]{rss}
C.~Rezk, S.~Schwede and B.~Shipley, 
{\em Simplicial structures on model categories and functors},
Amer.~J.~Math. {\bf 123} (2001), 551--575.

\bibitem[Ri03]{richter-DoldKan}
B.~Richter, {\em Symmetry properties of the Dold-Kan correspondence},
Math. Proc. Cambridge Philos. Soc. {\bf 134} (2003), 95--102.

\bibitem[Ric89a]{rickard1}
J.~Rickard, {\em Morita theory for derived categories}, 
J. London Math. Soc. (2) {\bf 39} (1989), 436--456.

\bibitem[Ric89b]{derived vs stable}
J.~Rickard, {\em Derived categories and stable equivalence}. 
J.~Pure Appl.\ Algebra {\bf 61} (1989), 303--317.

\bibitem[Ric91]{rickard2}
J.~Rickard, {\em Derived equivalences as derived functors}, 
J. London Math. Soc. (2) {\bf 43} (1991), 37--48.

\bibitem[Rob87a]{robinson-spectra}
A.~Robinson, {\em Spectral sheaves: a model category for 
stable homotopy theory}, J. Pure Appl. Algebra {\bf 45} (1987), 171--200.

\bibitem[Rob87b]{robinson-derived}
A.~Robinson, {\em The extraordinary derived category}, Math. Z. {\bf 196}
(1987), no.2, 231--238.

\bibitem[Sch97]{sch-cotangent}
S.~Schwede, {\em Spectra in model categories and applications to the 
algebraic  cotangent complex}, J. Pure Appl. Algebra {\bf 120} (1997), 77--104.

\bibitem[Sch01a]{sch-comparison}
S. Schwede, {\em {$S$}-modules and symmetric spectra}, Math.\ Ann.\ {\bf 319}
(2001), 517--532.

\bibitem[Sch01b]{sch-2}
S.~Schwede, {\em The stable homotopy category has a unique model at the
prime 2}, Adv.\ Math.\ {\bf 164} (2001), 24-40.

\bibitem[SS00]{ss-monoidal}
S.~Schwede and B.~Shipley, {\em Algebras and modules in monoidal 
model categories}, Proc.\ London Math.\ Soc.\ {\bf 80} (2000), 491-511.

\bibitem[SS02]{ss-uniqueness}
S.~Schwede and B.~Shipley, {\em A uniqueness theorem 
for stable homotopy theory}, Math.\ Z. {\bf 239} (2002), 803-828.

\bibitem[SS03]{ss-classification}
S.~Schwede and B.~Shipley, {\em Stable model categories are categories
of modules}, Topology {\bf 42} (2003), 103-153.

\bibitem[Sh2]{shipley-mon}
B.~Shipley, {\em Monoidal uniqueness of stable homotopy theory},
Adv. in Math. {\bf 160} (2001), 217-240.   

\bibitem[Spa88]{spaltenstein}
N.~Spaltenstein, {\em Resolutions of unbounded complexes},
Compositio Math.\ {\bf 65} (1988), 121--154.

\bibitem[TT]{Thomason-Trobaugh}
R.~W.~Thomason, T.~Trobaugh, 
{\em Higher algebraic $K$-theory of schemes and of derived categories},
The Grothendieck Festschrift, Vol. III, 247--435, 
Progr. Math., 88, Birkh{\"a}user Boston, Boston, MA, 1990.

\bibitem[Ver96]{verdier}
J.-L.~Verdier, {\em Des cat\'egories d\'eriv\'ees 
des cat\'egories ab\'eliennes}, {Ast\'erisque} {\bf 239} (1997).
With a preface by Luc Illusie, Edited and with a note by Georges Maltsiniotis,
xii+253 pp.

\bibitem[Voe]{voevodsky-milcon}
V.~Voevodsky, {\em The Milnor Conjecture}, Preprint (1996)

\bibitem[Voe98]{voevodsky-icm}
V.~Voevodsky, {\em ${\mathbb A}^1$-homotopy theory}, Doc. Math. ICM {\bf I},
1998, 417-442.
 
\bibitem[Wei94]{Weibel:1994a}
C.~A.~Weibel, {\em An introduction to homological algebra},
Cambridge Studies in Advanced Mathematics {\bf 38}. 
Cambridge University Press, Cambridge, 1994. xiv+450 pp.
\end{thebibliography}
\end{document}